\documentclass[number,citesort,dvips]{arxbj}

\aid{0}
\volume{17}
\issue{1}
\pubyear{2011}
\firstpage{88}
\lastpage{113}
\doi{10.3150/10-BEJ273}

\makeatletter
\newtheorem{theorem}{Theorem}[section]
\newtheorem{proposition}{Proposition}[section]
\newtheorem{corollary}{Corollary}[section]
\newtheorem{lemma}{Lemma}[section]

\newproclaim{definition}{Definition}[section]
\newproclaim{comment}{Comment}[section]

\newremark{remark}{Remark}[section]
\newremark{example}{Example}

\makeatother

\begin{document}
\begin{frontmatter}

\title{Invariance principles for linear processes with application to
isotonic regression}
\runtitle{Invariance principles for linear processes}

\begin{aug}
\author[a]{\fnms{J\'{e}r\^{o}me} \snm{Dedecker}\thanksref{a}\ead[label=e1]{jerome.dedecker@upmc.fr}},
\author[b]{\fnms{Florence} \snm{Merlev\`{e}de}\corref{}\thanksref{b}\ead[label=e2]{florence.merlevede@univ-mlv.fr}}\and
\author[c]{\fnms{Magda} \snm{Peligrad}\thanksref{c}\ead[label=e3]{peligrm@ucmail.uc.edu}}
\runauthor{J. Dedecker, F. Merlev\`{e}de and M. Peligrad}
\address[a]{Universit\'e Paris Descartes, Laboratoire MAP5, UMR CNRS 8145, 45 rue
 des Saints-P\`eres, F-75210 Paris cedex 06, France. \printead{e1}}
\address[b]{Universit\'e Paris Est-Marne la
Vall\'ee, LAMA and CNRS UMR 8050, 5 Boulevard Descartes, 77454 Marne
La Vall\'ee Cedex 2, France. \printead{e2}}
\address[c]{Department of Mathematical
Sciences, University of Cincinnati, P.O. Box 210025, Cincinnati, OH
45221-0025, USA. \printead{e3}}
\end{aug}

% HISTORY:
\received{\smonth{3} \syear{2009}}
\revised{\smonth{11} \syear{2009}}

% ABSTRACT
%
\begin{abstract}
In this paper, we prove maximal inequalities and study the functional central
limit theorem for the partial sums of linear processes generated by
dependent innovations. Due to the general weights, these processes can
exhibit long-range dependence and the limiting distribution is a fractional
Brownian motion. The proofs are based on new approximations by a linear
process with martingale difference innovations. The results are then applied
to study an estimator of the isotonic regression when the error process
is a
(possibly long-range dependent) time series.
\end{abstract}

% KEYWORDS
%
\begin{keyword}
\kwd{fractional Brownian motion}
\kwd{generalizations of martingales}
\kwd{invariance principles}
\kwd{isotonic regression}
\kwd{linear processes}
\kwd{moment inequalities}
\end{keyword}

\end{frontmatter}

%s1 ###
\section{Introduction and notation}

Without loss of generality, we assume that all the strictly stationary
sequences $(\xi_{i})_{i\in{\mathbf{Z}}}$ considered in this paper are
given by $\xi_{i}=\xi_{0}\circ T^{i}$, where $T\dvtx\Omega\mapsto
\Omega$ is
a bijective bimeasurable transformation preserving the probability
${\mathbf{%
P}}$ on $(\Omega,\mathcal{A})$. We denote by ${\mathcal{I}}$ the
$\sigma$%
-algebra of all $T$-invariant sets. For a subfield $\mathcal{F}_{0}$
satisfying $\mathcal{F}_{0}\subseteq T^{-1}(\mathcal{F}_{0})$, let ${%
\mathcal{F}}_{i}=T^{-i}({\mathcal{F}}_{0})$. Let $\mathcal{F}_{-\infty
}=\bigcap_{n\geq0}\mathcal{F}_{-n}$ and ${\mathcal{F}}_{\infty
}=\bigvee_{k\in{\mathbf{Z}}}{\mathcal{F}}_{k}$. The sequence
$({\mathcal{F}}%
_{i})_{i\in{\mathbf{Z}}}$ will be called a \textit{stationary filtration}.
We also assume that $\xi_{0}$ is \textit{regular}, that is, $\mathbf
{E}(\xi_{0}|{%
\mathcal{F}}_{-\infty})=0$ and $\xi_{0}$ is ${\mathcal{F}}_{\infty}$%
-measurable. On ${\mathbf{L}}^{2}$, we define the projection operator $P_{j}$
by
\begin{eqnarray*}
P_{j}(Y)=\mathbf{E}(Y|\mathcal{F}_{j})-\mathbf{E}(Y|\mathcal{F}_{j-1}).
\end{eqnarray*}
For any random variable $Y$, $\Vert Y\Vert_{p}$ denotes the norm in ${%
\mathbf{L}}^{p}$.

Recall that the linear process $X_{k}=\sum_{i\in{\mathbf{Z}}}a_{i}\xi
_{k-i} $ is well defined in ${\mathbf{L}}^{2}$ for any $(a_{i})_{i\in{%
\mathbf{Z}}}$ in $\ell^{2}$ (i.e., $\sum_{i\in{\mathbf
{Z}}}a_{i}^{2}<\infty$) if and only if the stationary sequence $(\xi
_{i})_{i\in{\mathbf{Z}}}$
has a bounded spectral density. Let $S_{n}=X_{1}+\cdots+X_{n}$ and $%
c_{n,j}=a_{1-j}+\cdots+a_{n-j}$. In the case where $\xi_{0}$ is
${\mathcal{%
F}}_{0}$-measurable, Peligrad and Utev \cite{PeligradUtev2006b}
have proven that if the
sequence $(\xi_{i})_{i\in{\mathbf{Z}}}$ satisfies an appropriate weak
dependence condition, then
\begin{eqnarray*}
\biggl(\sum_{j\in{\mathbf{Z}}}c_{n,j}^{2} \biggr)^{-1/2}S_{n}
\end{eqnarray*}
converges in distribution to $\sqrt{\eta}N$, where $\eta$ is a
non-negative ${%
\mathcal{I}}$-measurable random variable and $N$ is a standard normal
random variable independent of $\eta$. Their result extends the classical
result of Ibragimov \cite{Ibragimov1962} from i.i.d.~$\xi_{i}$'s
to the case of weakly
dependent sequences. In particular, the result applies if
%e1 ###
%
\begin{equation}\label{HeyHan}
\sum_{i\in{\mathbf{Z}}}\Vert P_{0}(\xi_{i})\Vert_{2}<\infty.
\end{equation}
Note that if this condition is satisfied, then the series $\sum_{k\in{%
\mathbf{Z}}}|{\mathbf{E}}(\xi_{0}\xi_{k})|$ converges. Indeed, since
$\xi_k = \sum_{i \in{%
\mathbf{Z}}} P_i ( \xi_k)$ and since ${\mathbf{E}} ( P_i ( \xi_0)
P_j(\xi_k))=0$ if $i \neq j$, it follows that for any $k \in
{\mathbf{Z}}$,
\[
|{\mathbf{E}}(\xi_{0}\xi_{k})| \leq \bigg|\sum_{i \in{%
\mathbf{Z}}} {\mathbf{E}} ( P_i ( \xi_0) P_i (\xi_k) ) \bigg| \leq\sum_{i
\in{%
\mathbf{Z}}} \Vert P_{0}(\xi_{i})\Vert_{2}\Vert P_{0}(\xi
_{k+i})\Vert_{2}
\]
so that $\sum_{k\in{%
\mathbf{Z}}}|{\mathbf{E}}(\xi_{0}\xi_{k})| \leq ( \sum_{i \in{%
\mathbf{Z}}} \Vert P_{0}(\xi_{i})\Vert_{2} )^2$. In addition,
under condition (\ref{HeyHan}), the non-negative random variable $\eta$
satisfies
$\eta= \sum_{k\in{%
\mathbf{Z}}}{\mathbf{E}}(\xi_{0}\xi_{k}| {\mathcal I})$.

Condition (\ref{HeyHan}) was introduced by Hannan \cite{Hannan1973}, and
by Heyde \cite{Heyde1974} in a slightly weaker form, and is well
adapted to
the analysis of time series (see, in particular, the application to
time series regression
given in the paper by Hannan \cite{Hannan1973}). As we shall see
in our Remark~\ref{use}%
, condition (\ref{HeyHan}) is also satisfied if
%e2 ###
%
\begin{equation}
\sum_{n=1}^{\infty}\frac{1}{\sqrt{n}}\Vert{\mathbf{E}}(\xi
_{n}|\mathcal{F}%
_{0})\Vert_{2}<\infty\quad\mbox{and}\quad\sum_{n=1}^{\infty}\frac{1}{\sqrt{n}}
\Vert\xi_{-n}-{\mathbf{E}}(\xi_{-n}|\mathcal{F}_{0})\Vert_{2}<\infty,
\label{racine}
\end{equation}
which is weaker than the condition introduced by Gordin \cite{Gordin1969}. If $\xi_{0}$
is ${\mathcal{F}}_{0}$-measurable, then condition (\ref{racine}) leads
to interesting new conditions for weakly dependent sequences and can be
successfully applied to functions of dynamical systems (see \cite{PeligradUtev2006b}, Section~3, and \cite{DMPU}, Section~6, for more details).

A natural question is now: what can we say about the weak convergence of
the partial sum process
%e3 ###
%
\begin{equation}\label{partialsumproc}
\biggl\{ \biggl(\sum_{j\in{\mathbf{Z}}}c_{n,j}^{2} \biggr)^{-1/2}S_{[nt]},t\in
[0,1] \biggr\}
\end{equation}
in the space $D([0,1])$ of cadlag functions equipped with the uniform
topology? Due to the results of Davydov \cite{Davydov1970} for
i.i.d.~$\xi_{i}$'s, we know
that the question is not as simple as for the central limit question and
that the limiting process (when it exists) depends on the behavior of the
normalizing sequence $v_{n}^{2}=\sum_{j\in{\mathbf{Z}}}c_{n,j}^{2}$. More
precisely, if (\ref{HeyHan}) holds and if there exists $\beta\in\,]0,2]$
such that
%e4 ###
%
\begin{equation}
\mbox{for any }t\,\in\,]0,1]\qquad \lim_{n\rightarrow\infty}\frac{%
v_{[nt]}^{2}}{v_{n}^{2}}=t^{\beta},
\end{equation}
then we show in Theorems \ref{IP} and \ref{IP2} that the finite-dimensional
marginals of the process (\ref{partialsumproc}) converge in distribution
to those of $\sqrt{\eta}W_{H}$, where $W_{H}$ is a fractional Brownian
motion, independent of $\eta$, with Hurst index $H=\beta/2$. The question
is now: under what conditions can we obtain the tightness in $D([0,1])$?

In Theorem \ref{IP} of Section \ref{sectIP2}, we show that if $\beta\in\,]1,2]$, then condition (\ref{HeyHan}) is sufficient for weak
convergence in $D([0,1])$. If $\beta\in\,]0,1]$, we point out in
Theorem \ref%
{sectIP2} that the convergence in $D([0,1])$ holds if (\ref{HeyHan}) is
replaced by the stronger condition
%e5 ###
%
\begin{equation}
\sum_{i\in{\mathbf{Z}}}\Vert P_{0}(\xi_{i})\Vert_{q}<\infty \qquad\mbox{for $q>2/\beta$}. \label{HeyHan0}
\end{equation}

As a matter of fact, for $\beta=1$, it is known from counterexamples given
in \cite{WuWoodroofe2004} and \cite{MerlevedePeligrad2006}
that if the sequence $(\xi_{i})_{i\in{\mathbf{Z}}}$ is i.i.d.~with ${%
\mathbf{E}}(\xi_{0}^{2})<\infty$, then the weak invariance principle may
not be true for the partial sums of the linear process, so a
reinforcement of (\ref{HeyHan}) is necessary. The case $\beta=1$,
where $%
W_{1/2}$ is a standard Brownian motion, is of special interest and is known
as the weakly dependent case. In that case, we point out in Section \ref
{sectIP3} that if we  impose some additional assumptions on
$(a_{i})_{i\in{%
\mathbf{Z}}}$, then condition (\ref{HeyHan}) is sufficient for the weak invariance
principle (Comments \ref{SM} and \ref{comH}) or may be reinforced in a
weaker way than (\ref{HeyHan0}) (Theorem \ref{IPlinearDMV}).

Note that, with the notation above, the sum $S_{n}$ may be written as
%e6 ###
%
\begin{equation}
S_{n}=\sum_{i\in{\mathbf{Z}}}c_{n,i}\xi_{i}. \label{deflinstat}
\end{equation}
Consequently, to prove our main theorems, in Section \ref{sectMom},
we give two preliminary results for linear statistics of type (\ref
{deflinstat}):
first, a moment inequality given in Proposition \ref{exp2} and, next, a
martingale approximation result given in Proposition \ref{exp3}, which
enables us to go back to the standard case where the $\xi_{i}$'s are
martingale differences. Both results are given in terms of Orlicz norms.

Our results provide, besides the invariance principles, estimates of the
maximums of partial sums that make them appealing to the study of statistics
involving linear processes. In Section \ref{sectISO}, we apply our
results to
the so-called isotonic regression problem
%e7 ###
%
\begin{equation}%\label{isotonicmodel}
y_{k}=\phi \biggl(\frac{k}{n} \biggr)+X_{k},\qquad k=1,2,\ldots,n,
\end{equation}
where $\phi$ is non-decreasing and the error $X_{k}$ is a linear
process. We follow the general scheme given in \cite{AnevskiHossjer2006}, who showed that in the
context of dependent
errors, the main tools to obtain the asymptotic distribution of the
isotonic estimator $\hat{\phi}$ are the
convergence in $D([0,1])$ of the partial sum process defined in (\ref{partialsumproc}) and a suitable maximal inequality for the
rescaled stochastic term (see their condition (14)). Zhao and
Woodroofe \cite{ZhaoWoodroofe2009} shed light on the fact that,
in addition to the weak
invariance principle, it is, in fact, enough to prove a suitable
maximal inequality directly on the partial sums of the error
process.
As in \cite{AnevskiHossjer2006}, the rate
of convergence of $\hat{\phi}$ is determined by
the asymptotic behavior of the normalizing sequence $v_{n}^{2}=\sum
_{j\in{%
\mathbf{Z}}}c_{n,j}^{2}$ and the limiting distribution depends on the
limiting process $W_{H}$.

%s2 ###
\section{Moment inequalities and martingale approximation for Orlicz norms}\label{sectMom}

For $\Psi\dvtx  {\mathbf{R}}_{+}\rightarrow{\mathbf{R}}_{+}$ a Young
function (convex, increasing, $\Psi(0)=0$ and $\lim_{x \rightarrow
\infty} \Psi(x) = \infty$), we denote by
${\mathbf{L}}_{\Psi}$ the Orlicz space defined as the space of all random
variables $X$ such that ${\mathbf{E}} \Psi(|X|/c)< \infty$ for some $c >0$.
It is a Banach space for the norm
\begin{eqnarray*}
\| X \|_{\Psi} = \inf\{ c >0, {\mathbf{E}} \Psi(|X|/c) \leq1 \}.
\end{eqnarray*}
Note that if $\Psi(x)=x^q$, $1 \leq q < \infty$, then ${\mathbf{L}}_{\Psi}={\mathbf{L}}^q$.

Let us also introduce the following class of functions (see \cite{DelapenaGine1999}, page~60). For $\alpha
>0$, the class $\mathcal{A}_{\alpha}$
consists of functions $\Phi\dvtx {\mathbf{R}}_{+}\rightarrow{\mathbf
{R}}_{+},$ where $%
\Phi(0)=0,$ $\Phi$ is non-decreasing continuous and such that
\begin{eqnarray*}
\Phi(cx)\leq c^{\alpha}\Phi(x)\qquad\mbox{for all }c\geq2, x\geq0.
\end{eqnarray*}
We also denote by ${\mathcal C} (\mathcal{A}_{\alpha})$ the class of
functions $\Psi$ such that $\Psi$ is a Young function in $\mathcal
{A}_{\alpha}$
and $x\mapsto
\Psi(\sqrt{x})$ is a convex function.
%p2.1
\begin{proposition}\label{exp2}
Let $\{Y_{k}\}_{k\in\mathbf{Z}}$ be a sequence of random
variables such that for all $k$, $\mathbf{E}(Y_{k}|{\mathcal
{F}}_{-\infty
})=0$ almost surely and $Y_{k}$ is ${\mathcal{F}}_{\infty
}$-measurable. Let
$\Psi$ be a function in ${\mathcal C} (\mathcal{A}_{\alpha})$. Assume that
%in $\mathcal{A}_{\alpha} $
%
\begin{eqnarray*}
\|P_{k-j}(Y_{k})\|_{\Psi}\leq p_{j}\quad \mbox{and}\quad D_{\Psi} :=\sum
_{j=-\infty}^{\infty}p_{j}<\infty.
\end{eqnarray*}
For any positive integer $m$, let $\{c_{m,j}\}_{j\in{\mathbf{Z}}}$ be a
sequence in $\ell^{2}$. Define $S_{m}=\sum_{j\in{\mathbf
{Z}}}c_{m,j}Y_{j}$%
. Then, for all $m\geq1$, there exists a positive constant $C_{\alpha}$,
depending only on $\alpha$, such that
%e8 ###
%
\begin{equation}%\label{r1lemmamom}
\Vert S_{m}\Vert_{\Psi}\leq C_{\alpha}D_{\Psi} \biggl(\sum_{j\in{\mathbf
{Z%
}}}c_{m,j}^{2} \biggr)^{1/2}.
\end{equation}
\end{proposition}
%r2.1
\begin{remark}
Using the notation of the above proposition, we get, for the special
function $%
\Psi(x)=x^{q}$ with $q\in[2,\infty[$, the following moment
inequality. Assume that
\begin{eqnarray*}
\|P_{k-j}(Y_{k})\|_{q}\leq p_{j} \quad\mbox{and}\quad D_{q}:=\sum_{j=-\infty}^{\infty}p_{j}<\infty.
\end{eqnarray*}
Then, for any $m\geq1$,
\begin{eqnarray*}
\|S_{m}\|_{q}\leq C_{q} \biggl(\sum_{j\in{\mathbf{Z}}}c_{m,j}^{2} \biggr)^{1/2}D_{q},
\end{eqnarray*}
where $C_{q}^{q}=18q^{3/2}/(q-1)^{1/2}$.
\end{remark}

For all $j\in{\mathbf{Z}}$, let $d_{j}=\sum_{\ell\in{\mathbf{Z}}%
}P_{j}(\xi_{\ell})$. Clearly, $(d_{j})_{j\in{\mathbf{Z}}}$ is a stationary
sequence of martingale differences with respect to the filtration
$(\mathcal{%
F}_{j})_{j\in{\mathbf{Z}}}$.
%p2.2
\begin{proposition}\label{exp3}For any positive integer $n$, let $\{c_{n,i}\}_{i\in{\mathbf{Z}}}$ be a sequence in $\ell^{2}$. Let
$\Psi$ be a function in ${\mathcal C} (\mathcal{A}_{\alpha})$. If
$\sum_{j\in\mathbf{Z}}\|P_{0}(\xi
_{j})\Vert_{\Psi}<\infty$, then we have the following
martingale-difference approximation: for any positive integer $m$,
there exists a
positive constant $C_{\alpha}$, depending only on $\alpha$, such that
\begin{eqnarray*}
\bigg\|\sum_{i\in{\mathbf{Z}}}c_{n,i}(\xi_{i}-d_{i}) \bigg\|_{\Psi} &\leq&2C_{\alpha} \biggl(\sum_{i\in{\mathbf{Z}}}c_{n,i}^{2} \biggr)%
^{1/2}\sum_{|k|\geq m}\Vert P_{0}(\xi_{k})\Vert_{\Psi} \\
&&{}+3C_{\alpha}m \biggl(\sum_{j\in{\mathbf{Z}}}(c_{n,j}-c_{n,j-1})^{2} \biggr)%
^{1/2}\sum_{j\in\mathbf{Z}}\|P_{0}(\xi_{j})\Vert_{\Psi}.
\end{eqnarray*}
%
%where $%
%d_{j}=\sum_{\ell\in{\mathbf{Z}}}P_{j}(\xi_{\ell})$ for all $j\in
%{\mathbf{Z}}$.
\end{proposition}

%c2.1
\begin{corollary}
\label{corexp3}Let $(a_{i})_{i\in{\mathbf{Z}}}$ be a sequence of
real numbers in $\ell^{2}$. Let $\Psi$ be a function in ${\mathcal
C} (\mathcal{A}_{\alpha})$. Assume that $\xi_{0}\in L_{\Psi}$ and
$\sum_{j}\|P_{0}(\xi_{j})\Vert_{\Psi}<\infty$.
Let $X_{k}=\sum_{j\in{\mathbf{Z}}}a_{j}\xi_{k-j}$ and $Y_{k}=\sum
_{j\in{%
\mathbf{Z}}}a_{j}d_{k-j}$. Set $S_{n}=\sum_{k=1}^{n}X_{k}$ and $%
T_{n}=\sum_{k=1}^{n}Y_{k}$. Then, for any positive $m$, there exist positive
constants $C_{1}$ and $C_{2}$ such that
%e9 ###
%
\begin{equation}
\Vert S_{n}-T_{n}\Vert_{\Psi}\leq C_{1}v_{n} \sum_{|k|\geq m}\Vert
P_{0}(\xi_{k})\Vert_{\Psi}+C_{2}m , \label{inecorexp3}
\end{equation}
where $v_{n}^{2}=\sum_{j\in{\mathbf{Z}}}c_{n,j}^{2}$ and $%
c_{n,j}=a_{1-j}+\cdots+a_{n-j}$.
\end{corollary}

\begin{pf} %{Proof of Corollary \ref{corexp3}}
We apply Proposition \ref{exp3} by noting that $S_{n}-T_{n}=\sum_{j\in{\mathbf{Z}}}c_{n,j}(\xi
_{j}-d_{j})$ and that
\begin{eqnarray*}
\sum_{j\in{\mathbf{Z}}}(c_{n,j}-c_{n,j-1})^{2}\leq4\sum_{j\in{\mathbf
{Z}}%
}a_{j}^{2}.
\end{eqnarray*}
\upqed
\end{pf}

Using the Orlicz norms, we give the following maximal inequality, which
is a
refinement of inequality (6) in \cite{Wu2007}, Proposition~1.
%l2.1
\begin{lemma}
\label{max} Let $\Psi$ be a Young function. Let $p\geq1$ and write
$\Psi_p(x)$
for $\Psi(x^p)$. Let $(Y_i)_{1 \leq i \leq2^N}$ be a strictly stationary
sequence of random variables such that $\|Y_1\|_{\Psi_p} < \infty$. Let $S_n
= Y_1 +\cdots+Y_n$. Then
\begin{eqnarray*}
\Big\| \max_{1 \leq m \leq2^N}|S_m| \Big\|_p \leq\sum_{L=0}^{N} \|
S_{2^L}\|_{\Psi_p} ( \Psi^{-1} ( 2^{N-L} ) )^{1/p}.
\end{eqnarray*}
\end{lemma}
%r2.2
\begin{remark}
Clearly, we can take $\Psi(x)=x$ in Lemma \ref{max}. Hence, in the
stationary case, we recover the inequality (6) in \cite{Wu2007}.
\end{remark}

%s3 ###
\section{Invariance principle for linear processes}\label{sectIP}

In this section, we shall focus on the weak invariance principle for linear
processes. Let $(a_{i})_{i\in{\mathbf{Z}}}$ be a sequence of real numbers
in $\ell^{2}$. Let
%e10 ###
%
\begin{equation} \label{def1}
X_{k}=\sum_{i\in{\mathbf{Z}}}a_{i}\xi_{k-i}\quad\mbox{and}\quad
S_{[nt]}=\sum_{k=1}^{[nt]}X_{k} ,
\end{equation}
and
%e11 ###
%
\begin{equation}\label{defvn}
v_{n}^{2}=\sum_{j\in{\mathbf{Z}}}c_{n,j}^{2},\qquad\mbox{where }
c_{n,j}=a_{1-j}+\cdots+a_{n-j}.
\end{equation}
The behavior of the process $\{S_{[nt]},t\in[0,1]\}$, properly
normalized, strongly depends on the behavior of the sequence
$(a_{i})_{i\in{%
\mathbf{Z}}}$.

In the next two sections, we treat separately the case where the limit
process is a mixture of fractional Brownian motions and the case where it
is a mixture of standard Brownian motions.

%s3.1 ###
\subsection{Convergence to a mixture of fractional Brownian motions}\label{sectIP2}

%d3.1
\begin{definition}
We say that a positive sequence $(v_{n}^{2})_{n\geq1}$ is regularly varying
with exponent $\beta>0$ if, for any $t\in\,]0,1]$,
%e12 ###
%
\begin{equation}\label{hyposn}
\frac{v_{[nt]}^{2}}{v_{n}^{2}}\rightarrow t^{\beta} \qquad\mbox{as }%
n\rightarrow\infty.
\end{equation}
\end{definition}

We shall separate the case $\beta\in\,]1,2]$ from the case $\beta\in\,]0,1].$
%t3.1
\begin{theorem}
\label{IP} Let $(a_{i})_{i\in{\mathbf{Z}}}$ be in ${\ell^{2}}$. Let
$\beta
\in\,]1,2]$ and assume that $v_{n}^{2}$ defined by (\ref{defvn}) is regularly
varying with exponent $\beta$. Let $\xi_{0}$ be a regular random variable
such that $\Vert\xi_{0}\Vert_{2}<\infty$ and let $\xi_{i}=\xi
_{0}\circ T^{i}.$ Assume that condition (\ref{HeyHan}) is satisfied.
The process $\{v_{n}^{-1}S_{[nt]},t\in[0,1]\}$ then converges
in $%
D([0,1]) $ to $\sqrt{\eta}W_{H}$, where $W_{H}$ is a standard fractional
Brownian motion independent of $\eta$ with Hurst index $H=\beta/2$, $%
\eta=\sum_{k\in{\mathbf{Z}}}{\mathbb{E}}(\xi_{0}\xi_{k}|{\mathcal{I}})$
and there exists a positive constant $C$ (not depending on $n$) such that
%e13 ###
%
\begin{equation}
{\mathbf{E}}\Bigl(\max_{1\leq k\leq n}S_{k}^{2}\Bigr)\leq Cv_{n}^{2}.
\label{inemawxIP}
\end{equation}
\end{theorem}
%t3.2
\begin{theorem}
\label{IP2} Let $\beta\in\,]0,1]$ and assume that $v_{n}^{2}$ defined
by (%
\ref{defvn}) is regularly varying with exponent $\beta$. Let $\xi_{0}$
be a
regular random variable such that $\Vert\xi_{0}\Vert_{2}<\infty$ and
let $\xi_{i}=\xi_{0}\circ T^{i}.$ Assume that condition (\ref
{HeyHan}) is
satisfied. The finite-dimensional distributions of $%
\{v_{n}^{-1}S_{[nt]},t\in[0,1]\}$ then converge to the
corresponding ones of $\sqrt{\eta}W_{H}$, where $W_{H}$ is a
standard fractional Brownian motion, independent of $\eta$, with
Hurst index $H=\beta/2$ and $\eta=\sum_{k\in
{\mathbf{Z}}}{\mathbf{E}}(\xi_{0}\xi_{k}|{\mathcal{I}})$. Assume,
in addition, that for a $q>2/{\beta}$, we have $\Vert\xi_{0}\Vert
_{q}<\infty$ and
%e14 ###
%
\begin{equation}
\sum_{j\in{\mathbf{Z}}}\Vert P_{0}(\xi_{j})\Vert_{q}<\infty.
\label{projcond1}
\end{equation}
Then the process $\{v_{n}^{-1}S_{[nt]},t\in[0,1]\}$ converges
in $%
D([0,1])$ to $\sqrt{\eta}W_{H}$ and (\ref{inemawxIP}) holds.
\end{theorem}
%r3.1
\begin{remark}
\label{var2}According to Peligrad and Utev \cite{PeligradUtev2006b}, Corollary 2, we have
\begin{eqnarray*}
\lim_{n\rightarrow\infty}\frac{\operatorname{Var}(S_{n})}{v_{n}^{2}}%
=\lim_{n\rightarrow\infty}\frac{\operatorname{Var}(\xi_{1}+\cdots+\xi
_{n})}{n}%
=v^{2}= \bigg\| \sum_{j \in{\mathbf{Z}}}P_{0}(\xi_{j}) \bigg\| _{2}^{2}.
\end{eqnarray*}
\end{remark}
%r3.2
\begin{remark}
\label{necess}In the context of Theorem \ref{IP}, condition (\ref{hyposn})
is necessary for the conclusion of this theorem (see \cite{Lamperti1962}). This
condition has also been imposed by Davydov \cite{Davydov1970} to
study the weak
invariance principle of linear processes with i.i.d.~innovations. To be
more precise, Davydov proved
that if (\ref{hyposn}) holds and if $\xi_0 \in{\mathbb L}^q$ with $q
\geq4$
and $q > 4 (1/\beta-1)$, then $\{v_{n}^{-1}S_{[nt]},t\in[0,1]\}
$ converges in $%
D([0,1])$ to $\sqrt{{\mathbf E}(\xi_0^2)}W_{\beta/2}$. Later, in the
case $\beta>1$, Konstantopoulos and Sakhanenko \cite{KonstantopoulosSakhanenko2004} sharpened
Davydov's result, showing that the weak invariance principle holds if
the $\xi_i$'s are i.i.d.~and in ${\mathbf L}^2$.
\end{remark}
%ex1
\begin{example}
For $0<d<1/2$, let us consider the linear
process $X_{k}$
defined by
%e15 ###
%
\begin{equation}
X_{k}=(1-B)^{-d}\xi_{k}=\sum_{i\geq0}a_{i}\xi_{k-i} , \label{deffractlin}
\end{equation}
where $B$ is the lag operator, $a_0 = 1$, $a_{i}=%
\frac{\Gamma(i+d)}{\Gamma(d)\Gamma(i+1)}$ for $i\geq1$ and $(\xi
_{i})_{i\in{\mathbf{Z}}%
} $ is a strictly stationary sequence satisfying the condition of
Theorem %
\ref{IP}. In this case, Theorem \ref{IP} applies with $\beta=2d+1$
since $%
a_{k}\sim(\Gamma(d))^{-1}k^{d-1}$.
\end{example}

%ex2
\begin{example}
Now, consider the following choice of
$(a_{k})_{k\geq0}$: $a_{0}=1$ and $a_{i}=(i+1)^{-\alpha}-i^{-\alpha
}$ for
$i\geq1$ with $\alpha\in\,]0,1/2[.$ Theorem \ref{IP2} then applies. Indeed,
for this choice, $v_{n}^{2}\sim\kappa_{\alpha}n^{1-2\alpha}$, where
$%
\kappa_{\alpha}$ is a positive constant depending on $\alpha$.
\end{example}
%ex3
\begin{example}
For the choice $a_{i}\sim i^{-\alpha}\ell
(i)$, where $\ell$ is a slowly varying function at infinity and
$1/2<\alpha
<1$, we have $v_{n}^{2}\sim\kappa_{\alpha}n^{3-2\alpha}\ell
^{2}(n)$ (see,
e.g., \cite{WangLinGulati2001}, relations (12)), where
$\kappa
_{\alpha}$ is a positive constant depending on $\alpha$.
\end{example}
%ex4
\begin{example}
Finally, if $a_{i}\sim i^{-1/2}(\log
i)^{-\alpha}$ for some $\alpha>1/2$, then $v_{n}^{2}\sim n^{2}(\log
n)^{1-2\alpha}/(2\alpha-1)$ (see \cite{WangLinGulati2001}, relations (12)). Hence, (\ref{hyposn}) is satisfied with
$\beta=2$.
\end{example}

For the sake of applications, we now give a sufficient condition
for (\ref{projcond1}) to hold.
%r3.3
\begin{remark}
\label{use}For any $q\in[2,\infty[$, the condition (\ref
{projcond1}) is satisfied if we assume that
%e16 ###
%
\begin{equation} \label{Mmax}
\sum_{n=1}^{\infty}\frac{1}{n^{1/q}}\Vert{\mathbf{E}}(\xi
_{n}|\mathcal{F}%
_{0})\Vert_{q}<\infty\quad\mbox{and}\quad\sum_{n=1}^{\infty}\frac{1}{n^{1/q}}
\Vert\xi_{-n}-{\mathbf{E}}(\xi_{-n}|\mathcal{F}_{0})\Vert_{q}<\infty.
\end{equation}
\end{remark}

The fact that (\ref{Mmax}) implies (\ref{projcond1}) extends \cite{PeligradUtev2006b}, Corollary 2, and also
%Dedecker, Merlev\`ede and Voln\'y
\cite{DMPU}, Corollary 5, from the case $q=2$ to more
general situations.

For causal linear processes, Shao and Wu \cite{ShaoWu2006} also
showed that the weak
invariance principle holds under the condition (\ref{projcond1}), as
long as
the coefficients of the linear processes satisfy a certain regularity
condition. To be more precise, their condition on the coefficients of the
linear processes lead either to $\beta>1$ or $\beta<1$. For this last
case, they specified the coefficients $(a_{i})_{i\geq0}$ as follows:
for $%
1<\alpha<3/2$, $a_{j}=j^{-\alpha}\ell(j)$ for $j\geq1$ (where $\ell(i)$
is a slowly varying function) and $\sum_{j=0}^{\infty}a_{j}=0$ (see,
e.g., their Lemma 4.1). For this choice, $v_{n}^{2}$ is regularly
varying with coefficient $\beta=3-2\alpha<1$. Our Theorem \ref{IP2} does
not require conditions on the coefficients, but only the fact that the
variance is regularly varying, which is a necessary condition.

%s3.2 ###
\subsection{Convergence to a mixture of Brownian motions\label{sectIP3}}

The case $\beta=1$ deserves special attention. For this case, the
limit is a
mixture of Brownian motions.

As an immediate consequence of Theorem \ref{IP2}, we formulate the following
corollary for causal linear processes, under a recent condition introduced
by Wu and Woodroofe \cite{WuWoodroofe2004}.
%c3.1
\begin{corollary}
\label{IPWW} Let $\xi_{0}$ be a regular random variable such that
$\Vert\xi_{0}\Vert_{q}<\infty$ for some $q >2$ and let $\xi_{i}=\xi_{0}\circ T^{i}$. Assume, in addition, that
%e17 ###
%
\begin{equation}
\sum_{j\in{\mathbf{Z}}}\Vert P_{0}(\xi_{j})\Vert_{q}<\infty.
\label{cond2delta}
\end{equation}
Let $(a_{i})_{i\in{\mathbf{Z}}}$ be a sequence of real numbers in $\ell
^{2} $ such that $a_{i}=0$ for $i<0$. Let $b_{j}=a_{0}+\cdots+a_{j}$. Define
$(X_{k})_{k\geq1}$ as above and assume that
%e18 ###
%
\begin{equation}
\sum_{k=0}^{n-1}b_{k}^{2}\rightarrow\infty\qquad\mbox{as }n\rightarrow
\infty,
\label{condbjinfty}
\end{equation}
and that
%e19 ###
%
\begin{equation}
\sum_{j=0}^{\infty}(b_{n+j}-b_{j})^{2}=\mathrm{o} \Biggl(\sum_{k=0}^{n-1}b_{k}^{2}%
\Biggr). \label{condbj}
\end{equation}
Then $v_{n}^{2}\sim nh(n)$, where $h(n)$ is a slowly varying function.
Moreover, the process $\{v_{n}^{-1}S_{[nt]},t\in[0,1]\}$ converges
in $D([0,1])$ to $\sqrt{\eta}W$, where $W$ is a standard Brownian motion,
independent of~$\eta$, and $\eta=\sum_{k\in{\mathbf{Z}}}{\mathbf
{E}}(\xi
_{0}\xi_{k}|{\mathcal{I}})$. In addition, (\ref{inemawxIP}) holds.
\end{corollary}

To prove this result, it suffices to apply Theorem \ref{IP2} and to use the
fact that under (\ref{condbjinfty}) and (\ref{condbj}), $v_{n}^{2}\sim nh(n)$
(see \cite{WuWoodroofe2004}). Under the same
conditions (\ref{condbjinfty})
and (\ref{condbj}), Wu and Min \cite{WuMin2005}, in their Theorem
1, also proved the weak
invariance principle, but under the stronger condition $\sum_{j\geq
0}j\Vert
P_{0}(\xi_{j})\Vert_{q}<\infty$ (in their paper, the random
variables $\xi
_{j}$ are adapted to the filtration ${\mathcal{F}}_{j}$).
%r3.4
\begin{remark}
The above result fails if, in (\ref{cond2delta}), we take $q=2$;
see \cite{WuWoodroofe2004} and also \cite{MerlevedePeligrad2006},
Example 1, page~657.
\end{remark}

Let us make some comments on the case where the condition (\ref{HeyHan}) is
sufficient for weak convergence to the Brownian motion with the
normalization $\sqrt n$. The first case is already known and the second
case deserves a short proof.
%c3.1
\begin{comment}
\label{SM} When $\sum_{i\in{\mathbf{Z}}}|a_{i}|<\infty$ (the
short memory case) and condition (\ref{HeyHan}) is satisfied, one can
use the result from \cite{PeligradUtev2006a} in
the adapted case,
showing that the invariance principle for the linear process is
inherited from the innovations at no extra cost. For this case, the
process $\{n^{-1/2}S_{[nt]},t\in[0,1]\}$ converges in
distribution in $D([0,1])$ to $\sqrt{\eta}W$, where $W$ is a
standard Brownian motion, independent of ${\eta}$, and $\eta
=A^{2}\sum_{k\in{\mathbf{Z}}}{\mathbf{E}}(\xi_{0}\xi
_{k}|{\mathcal{I}})$ with $A=\sum_{i\in{\mathbf{Z}}}a_{i}$.
Moreover, ${\mathbf{E}}(\max_{1\leq k\leq n}S_{k}^{2})\leq Cn.$ See \cite{DMPU},
Corollaries 2 and 3, for the non-adapted case.
\end{comment}
%c3.2
\begin{comment}
\label{comH} Let $(a_{i})_{i\in Z}$ in $\ell^2$ and assume that the
series $%
\sum_{i\in{\mathbf{Z}}}a_{i}$ converges (meaning that the two series $%
\sum_{i\geq0}a_{i}$ and $\sum_{i<0}a_{i}$ converge) and Heyde's
\cite{Heyde1975}
condition \textup{(H)} holds:
%
\begin{eqnarray*}
\mathrm{(H)}\quad\sum_{n=1}^{\infty} \biggl(\sum_{k\geq n}a_{k} \biggr)^{2}<\infty
\quad\mbox{and}\quad \sum_{n=1}^{\infty}
\biggl(\sum_{k\leq-n}a_{k} \biggr)^{2}<\infty.
\end{eqnarray*}
%
Assume, also, that condition (\ref{HeyHan}) is satisfied. The same
conclusion as in Comment \ref{SM} then holds.
\end{comment}

%ex5
\begin{example}
Heyde's condition allows the following
possibility: $\sum_{i\in{\mathbf{Z}}}|a_{i}|=\infty$, but $\sum_{i\in{
\mathbf{Z}}}a_{i}$ converges. For instance, if, for $n<0$, $a_{n}=0$
and, for
$n\geq1$, $a_{n}=(-1)^{n}u_{n}$, for some sequence $(u_{n})_{n\geq1}$ of
positive coefficients decreasing to zero such that $\sum_{n\geq
1}u_{n}=\infty$, then condition \textup{(H)} is satisfied as soon as $%
\sum_{n>0}u_{n}^{2}<\infty$, which is a minimal condition. It is noteworthy
to indicate that Heyde's condition implies (\ref{condbj}).
\end{example}

Now, if $\sum_{j\in{\mathbf{Z}}}|a_{j}|=\infty$ and \textup{(H)} does not hold,
then condition (\ref{cond2delta}) may still be weakened in some
particular cases.
The following result generalizes Corollary 4 in \cite{DMPU} to the case where the innovations of the
linear process are
not necessarily martingale difference sequences. We write
%e20 ###
%
\begin{equation}
s_{n}^{2}=n \Biggl(\sum_{i=-n}^{n}a_{i} \Biggr)^{2}. \label{defsn}
\end{equation}

%t3.3
\begin{theorem}
\label{IPlinearDMV} Let $(a_{i})_{i\in{\mathbf{Z}}}$ be a sequence of real
numbers in $\ell^{2}$, but not in $\ell^{1}$, and let $s_{n}^{2}$ be
defined by (\ref{defsn}). Define $(X_{k})_{k\geq1}$ as above and assume
that
%e21 ###
%
\begin{equation}\label{condDMV}
\limsup_{n\rightarrow\infty}\frac{\sum_{i=-n}^{n}|a_{i}|}{ |%
\sum_{i=-n}^{n}a_{i} |}<\infty\quad\mbox{and}\quad \sum_{k=1}^{n}\sqrt{%
\sum_{|i|\geq k}a_{i}^{2}}=\mathrm{o}(s_{n}).
\end{equation}
If one of the following two conditions holds,\vspace*{6pt}

(\textup{a})\quad $ \displaystyle \sum_{j\in{\mathbf{Z}}}\Vert P_{0}(\xi
_{j})\Vert_{\Psi_{2,\alpha}}<\infty$,\qquad
where $\Psi_{2,\alpha}(x)=x^{2}\log^{\alpha}(1+x^{2})$ and $\alpha>2$,\vspace*{6pt}

\noindent or\vspace*{6pt}

(\textup{b})\quad $ \displaystyle \sum_{j\in{\mathbf{Z}}}\log(1+|j|)\Vert P_{0}(\xi
_{j})\Vert_{2}<\infty$,\vspace*{6pt}

\noindent then $\{s_{n}^{-1}S_{[nt]},t\in[0,1]\}$ converges weakly in $%
D([0,1]) $ to $\sqrt{\eta}W$, where $W$ is a standard Brownian motion,
independent of $\eta$, and $\eta=\sum_{k\in{\mathbf{Z}}}{\mathbf
{E}}(\xi
_{0}\xi_{k}|{\mathcal{I}})$. In addition, there exists a positive
constant $%
C$ (not depending on $n$) such that
%e22 ###
%
\begin{equation}\label{inemawxIPlinearDMV}
{\mathbf{E}}\Bigl(\max_{1\leq k\leq n}S_{k}^{2}\Bigr)\leq Cs_{n}^{2}.
\end{equation}
\end{theorem}
%r3.5
\begin{remark}
\label{rmsn} For two positive sequences of numbers, the notation
$u_{n}\sim
v_{n}$ means that $\lim_{n\rightarrow\infty}u_{n}/v_{n}=1$. According to
\cite{DMPU}, Remark 12, we
have that
\begin{eqnarray*}
s_{n}^{2}\sim v_{n}^{2}\sim nh(n) ,
\end{eqnarray*}
where $h(n)$ is a slowly varying function at infinity. In addition, if we
assume the first part of condition (\ref{condDMV}) and $\sum_{j\in
{\mathbf{Z%
}}}|a_{j}|=\infty$, then we get that $s_{n}/\sqrt{n}\rightarrow\infty
$ as $%
n\rightarrow\infty$.
\end{remark}

%ex6
\begin{example}\label{ex6}
Consider the following choice of $%
(a_{k})_{k\in{\mathbf{Z}}}$: $a_{0}=1$ and $a_{i}=1/|i|$ for $i\neq0$.
Then Theorem \ref{IPlinearDMV} applies. Indeed, for this choice, condition
(\ref{condDMV}) holds and $s_{n}\sim2\sqrt{n}(\log n)$.
\end{example}

We now give a useful sufficient condition for the validity of condition
(b) of Theorem \ref{IPlinearDMV}.

%r3.6
\begin{remark}
\label{corIPlinearDMV} Condition (b) of Theorem \ref{IPlinearDMV} is
satisfied if we assume that
%e23 ###
%
\begin{equation} \label{Mmaxlog}
\sum_{n=1}^{\infty} \log n\frac{\Vert{\mathbf{E}}(\xi_{n}|\mathcal{F}_{0})\Vert_{2}}{\sqrt n}<\infty
\quad\mbox{and}\quad
\sum_{n=1}^{\infty} \log n \frac{\Vert\xi_{-n}-{\mathbf{E}}(\xi_{-n}|\mathcal{F}_{0})\Vert_{2}}{\sqrt n}<\infty.
\end{equation}
\end{remark}

%s4 ###
\section{Application to isotonic regression}\label{sectISO}

Let $\phi$ be a non-decreasing function on the unit interval and let
%e24 ###
%
\begin{equation} \label{isotonicmodel}
y_k=\phi \biggl( \frac kn \biggr) + X_k,\qquad k =1,2,\ldots,n ,
\end{equation}
where $(X_k)$ is a strictly stationary sequence of random variables
such that ${\mathbf{E}}(X_k)=0$ and ${\mathbf{E}}(X_k^2)<\infty$.
The problem is then to estimate $\phi$ in a nonparametric way. We write
$S_n = \sum_{k=1}^n X_k$.

Taking advantage of the monotonicity of the regression function,
isotonic estimates have been suggested. Let $\mu_{k}=\phi(k/n)$.
It is well known that the least-squares estimator
\begin{eqnarray*}
\hat{\mu}=\mathop{\mathrm{argmin}} \Biggl\{\sum_{k=1}^{n}(y_{k}-\mu_{k})^{2},\mu
_{1}\leq\cdots\leq\mu_{n} \Biggr\}
\end{eqnarray*}
is such that
\begin{eqnarray*}
\hat{\mu}_{k}=\max_{i\leq k}\min_{j\geq k}\frac{y_{i}+\cdots+y_{j}}{j-i+1}.
\end{eqnarray*}
In addition, setting
\begin{eqnarray*}
Y_{n}(t)=\frac{1}{n} \Biggl(\sum_{k=1}^{[nt]}y_{k} \Biggr)\quad\mbox{and}\quad\widetilde{Y%
}_{n}=\mbox{GCM}(Y_{n}) ,
\end{eqnarray*}
where GCM designates the greatest convex minorant, we have
\begin{eqnarray*}
\hat{\mu}_{k}=\widetilde{Y}_{n}^{\prime} \biggl(\frac{k}{n} \biggr) ,
\end{eqnarray*}
where the derivative in taken on the left (see
\cite{RobertsonWrightDykstra1988}). Now, let $\hat{\phi}_{n}(\cdot)$
be the left-continuous step function on
$[0,1]$ such that $\hat{\phi}_{n}(k/n)=\hat{\mu}_{k}$ at the knots
$k/n$ for
$k=1,\ldots,n$.

When the error process $(X_k)$ in the model (\ref{isotonicmodel}) is
short-range dependent and satisfies suitable weak dependence
conditions, Zhao and Woodroofe \cite{ZhaoWoodroofe2009} have
obtained the asymptotic behavior
of $\hat\phi_n (t)$. In their paper, an application to global
warming is given. Some other situations are considered in
\cite{AnevskiHossjer2006}: in their Theorem 3(iii),
they consider the case where $(X_k)$ can exhibit long-range
dependence and they assume that $X_k$ is a function
of a Gaussian process such that its Hermite polynomial expansion is
of rank greater than one. When no shape assumption is imposed on the regression
function, nonparametric regression analysis when data can exhibit long-range
dependence has been also studied by other authors (see, e.g.,
\cite{Robinson1997} or, more
recently, Gao and Wang \cite{GaoWang2006} wherein random
designs are introduced in the nonparametric trend model). The
motivation for
studying such models is that, in order to avoid misrepresenting the
mean function or the conditional mean function of long-range
dependent data, one should let the data ``speak for themselves'' in
terms of specifying the true form of the mean function or the
conditional mean function. Situations where the error process
$(X_k)$ in the model (\ref{isotonicmodel}) is long-range dependent
often occur when considering financial or climatological time series.
For instance, the annual series of winter means of the NAO index
(North Atlantic Oscillation index) exhibits long-range dependence
(see \cite{StephensonPavanBojariu2000}) and also an increasing trend
for the last decade (which can possibly be explained by global
warming). Concerning financial time series, we refer to the paper by
Pesee \cite{Pesee2008}, where daily exchange rate data are
studied. For
instance, the daily changes of the US dollar against the Deutsche Mark
constitute a financial series that exhibits long-range dependence with
a long
period of monotonic trend.
For other data examples of long-memory processes, we refer to the book
by Beran \cite{Beran1994}.\vadjust{\goodbreak}
In particular, concerning the monthly temperature for the northern hemisphere,
Beran suggests (page~29 of his book) that the series could be
long-range dependent
(see his Figure 1.12a--c, page~31).

The aim of this section, then, is to derive the asymptotic behavior of
$\hat\phi_n (t)$ when $X_k$ is a linear process which can exhibit
short or long memory. Recall that, by the well-known Wold
decomposition, a stationary process in ${\mathbf L}^2$ that is
purely non-deterministic and such that its one-step mean squared
error is positive can be represented by a linear process generated
by orthogonal random variables.

As is implicitly mentioned in \cite{AnevskiHossjer2006} and
elucidated in \cite{ZhaoWoodroofe2009}, the
two main tools to
obtain the asymptotic behavior of $\hat\phi_n (t)$ are a weak
invariance principle for the partial sums process $\{ S_{[nt]} , t
\in[0,1] \}$, properly normalized, and a suitable moment inequality
for $\max_{1 \leq k \leq n} S^2_k$.

%t4.1
\begin{theorem}
\label{thmiso1} Let $(a_{i})_{i\in{\mathbf{Z}}}$ and $(\xi_{i})_{i\in
{%
\mathbf{Z}}}$ be as in Comments \ref{SM} or \ref{comH}. Let us consider the
model~(\ref{isotonicmodel}) with $X_{k}$ defined by (\ref{def1}). For
any $%
t\in(0,1)$ such that $\phi^{\prime}(t)>0$,
\begin{eqnarray*}
n^{1/3}\kappa^{-1}\bigl(\hat{\phi}_{n}(t)-\phi(t)\bigr)\quad\Longrightarrow\quad\bigl(\sqrt{\eta}\bigr)^{2/3}\mathop{\mathrm{argmin}}\{B(s)+s^{2},s\in{\mathbf{R}}\} ,
\end{eqnarray*}
where $B$ denotes a standard two-sided Brownian motion independent of
$\eta$%
, $\eta=\sum_{k\in{\mathbf{Z}}}{\mathbf{E}}(\xi_{0}\xi_{k}|{\mathcal
{I}}%
) $ and $\kappa=2 (\frac{1}{2}A^{2}\phi^{\prime}(t) )^{1/3}$
with $A=\sum_{j\in{\mathbf{Z}}}a_{j}$.
\end{theorem}

%Let $v_n^2$ be defined by (\ref{defvn}). Assume that $v_n^2$ is
%regularly varying at infinity with exponent $\beta\in\,]0, 2]$; that
%is
%v_n^2 = n^{\beta} h(n) \mbox{ where $h(x)$ is a slowly varying
%function at infinity}.
Let $\beta\in\,]0, 2]$ and let $h$ be a slowly varying function at
%e25 ###
infinity. Now, let
\begin{equation} \label{defL}
L(x) = \biggl( \frac{1}{h (x^{2/(4 - \beta)})} \biggr)^{1/2}
\end{equation}
and note that $L(x)$ is also a slowly varying function at infinity.
Denote by $L^{*}$ the asymptotic conjugate of $L$, which means that
$L^{*}$ satisfies
%e26 ###
%
\begin{equation} \label{defL*}
\lim_{x \rightarrow\infty} L^{*} (x) L (x L^* (x) ) = 1.
\end{equation}
%
%Then as soon as we get the representation (\ref{hypovnhn}) for
%$v_n^2$,
Then define
%e27 ###
\begin{equation} \label{defdn}
d_n = \frac{1}{n^{(2-\beta)/(4-\beta)}} \ell(n), \qquad\mbox{where }\ell(n) = ( L^* (n) )^{2/(4 - \beta)}.
\end{equation}
%t4.2
\begin{theorem}
\label{thmiso2} Let $(a_i)_{i \in{\mathbf{Z}}}$ and $(\xi_{i})_{i \in
{%
\mathbf{Z}}} $ be as in Theorem \ref{IPlinearDMV}. For $\beta=1$ and $h(n)
=| \sum_{i=-n}^n a_i |^2$, let $d_n$ be defined by (\ref{defdn}). Let us
consider the model (\ref{isotonicmodel}) with $X_k$ defined by (\ref{def1}).
For any $t \in(0,1)$ such that $\phi^{\prime}(t) >0$,
\begin{eqnarray*}
d_n^{-1} \kappa^{-1} \bigl(\hat\phi_n (t) - \phi(t) \bigr)\quad \Longrightarrow\quad
\bigl(\sqrt{\eta} \bigr)^{2/3}\mathop{\mathrm{argmin}} \{B(s) + s^2 , s \in
{\mathbf{R}} \} ,
\end{eqnarray*}
where $B$ denotes a standard two-sided Brownian motion independent of
$\eta$%
, $\eta= \sum_{k \in{\mathbf{Z}}} {\mathbf{E}} (\xi_0 \xi_k |{\mathcal
{I}}%
) $ and $\kappa= 2 ( \frac12 \phi^{\prime}(t) )^{1/3} $.
\end{theorem}
%ex7
\begin{example}
In the case of the linear process defined in Example~\ref{ex6}, Theorem~\ref{thmiso2} applies with $d_n = n^{-1/3} ( 4 \ln(n) /3 )^{2/3}$.
\end{example}
%t4.3
\begin{theorem}
\label{thmiso3} Let $(a_i)_{i \in{\mathbf{Z}}}$ and $(\xi_{i})_{i \in{%
\mathbf{Z}}} $ be as in Theorem \ref{IP} or \ref{IP2} for some $\beta
\in\,]0, 2[$. By assumption, $v_n^2 $ defined by (\ref{defvn}) is regularly
varying with exponent $\beta$. For this $\beta$ and for $h(n)=v_n^2
n^{-\beta}$, let $d_n$ be defined by (\ref{defdn}). Let us consider the
model (\ref{isotonicmodel}) with $X_k$ defined by (\ref{def1}). Then,
for any
$t \in(0,1)$ such that $\phi^{\prime}(t) >0$,
we have
\begin{eqnarray*}
d_n^{-1} \kappa_{\beta}^{-1} \bigl(\hat\phi_n (t) - \phi(t) \bigr)\quad
\Longrightarrow\quad
\bigl(\sqrt{\eta} \bigr)^{1/(2-H)} \mathop{\mathrm{argmin}} \{B_H(s) + s^2, s \in{\mathbf{R}}
\} ,
\end{eqnarray*}
where $B_H$ denotes a standard two-sided fractional Brownian motion,
independent of $\eta$, with Hurst index $H= \beta/2$, $\eta= \sum_{k
\in{%
\mathbf{Z}}} {\mathbf{E}} (\xi_0 \xi_k |{\mathcal{I}})$ and where the
constant $%
\kappa_\beta$ is given by $\kappa_{\beta} =2 ( \phi^{\prime}(t)/2
)^{(2-\beta)/(4-\beta)} $.
\end{theorem}

%ex8
\begin{example}
In the case of the linear process defined
in Example 1, Theorem \ref{thmiso3} applies with $\beta= 2d +1$ and
$d_n = \tau_d n^{(1-2d)/(3-2d)}$, where $\tau_d$ is a positive
constant depending only on $d$.
\end{example}

\begin{pf*}{Proofs of Theorems \ref{thmiso1}--\ref{thmiso3}} For any $t\in(0,1)$ and any $s\in[
-td_{n}^{-1},d_{n}^{-1}(1-t)]$, let
\begin{eqnarray*}
Z_{n}(s)=d_{n}^{-2} \bigl(Y_{n}(t+d_{n}s)-Y_{n}(t)-\phi(t)d_{n}s \bigr).
\end{eqnarray*}
Then $d_{n}^{-1}(\hat{\phi}_{n}(t)-\phi(t))=\widetilde{Z}_{n}^{\prime
}(0)$%
, the left-hand derivative of the GCM of $Z_{n}$ at $s=0$. Hence, the
key for
establishing the result is the study of the GCM of the process $Z_{n}$. This
can be done by following the arguments given in
\cite{AnevskiHossjer2006}, Section 3, and also in \cite{ZhaoWoodroofe2009}. More precisely, a careful
analysis of the proofs given in
both of these papers shows that the following lemma is valid.

%l4.1
\begin{lemma} \label{AH}
Assume that there exists a positive sequence $m_{n}\rightarrow\infty$
satisfying, for any $t>0$,
%e28 ###
%
\begin{equation}
m_{[nt]}/m_{n}\rightarrow t^{H},\qquad\mbox{ where $H\in\,]0,1[$}, \label{condmn}
\end{equation}
and such that:
\begin{enumerate}[(1)]
\item[(1)] the process $\{m_{n}^{-1}S_{[nt]},t\in[0,1]\}$ converges
in $%
D([0,1])$ to $\sqrt{\eta}W_{H}$, where $\eta$ is a positive random
variable and $W_{H}$ is a standard fractional Brownian motion (with Hurst
index $H$) independent of $\eta$;

\item[(2)]${\mathbf E}(\max_{1\leq k\leq n}S_{k}^{2})\leq Cm_{n}^{2}$.
\end{enumerate}
Then, for any positive sequence $d_{n}\rightarrow0$ such that $%
nd_{n}\rightarrow\infty$ and $d_{n}^{-2}n^{-1}m_{[nd_{n}]}\rightarrow1$,
and for any $t\in(0,1)$ such that $\phi^{\prime}(t)>0$,
\begin{eqnarray*}
d_{n}^{-1}\kappa_{H}^{-1}\bigl(\hat{\phi}_{n}(t)-\phi(t)\bigr)\quad\Longrightarrow\quad\bigl(\sqrt
{%
\eta}\bigr)^{1/(2-H)}\mathop{\mathrm{argmin}}\{B_{H}(s)+s^{2},s\in
{\mathbf{R}}\} ,
\end{eqnarray*}
where $B_{H}(\cdot)$ denotes a standard two-sided fractional Brownian motion,
independent of ${\eta}$, with Hurst index $H\in\,]0,1[$ and $\kappa_{H}=2(\phi^{\prime}(t)/2 )^{(1-H)/(2-H)}$.
\end{lemma}

\begin{pf}%{Proof of Lemma \ref{AH}}
We proceed as in the proof of Anevski and H\"{o}ssjer \cite{AnevskiHossjer2006}, Theorem
3. The main point is then to
verify their assumptions A1--A7 in order to apply their Corollary 1.
Since $n d_n \rightarrow\infty$, assumption A2 follows from the
arguments given in the proof of Anevski and
H\"{o}ssjer \cite{AnevskiHossjer2006}, Theorem 3(i). By the
properties of our limiting process,
$\sqrt{\eta}W_{H}$, the assumptions A5 and A7 are satisfied. Now, if
assumption~A1 holds, then, by Anevski and H\"{o}ssjer \cite{AnevskiHossjer2006}, Proposition~2,
and the properties of the fractional Brownian motion, assumption
A6 also holds. Note that their Proposition 2 allows the
continuous mapping theorem to be applied to the functional $h$ from
$D[-c,c]$ (the
space of cadlag functions on $[-c,c]$) to $\mathbb R$, defined as the
left-hand derivative of $\mbox{GCM$(x)$}$ at $0$. To verify their
assumptions A3 and A4, it suffices to apply their Proposition 1.
According to the proofs of their Lemmas B1 and B2, the condition
(14) of their Proposition 1 is satisfied as soon as their condition
(87) and our condition
(\ref{condmn}) are. Now, their condition (87) is clearly satisfied
provided item 2 of Lemma \ref{AH} holds.

It remains to prove \cite{AnevskiHossjer2006}, assumption A1, namely, that the process
\[
\bigl\{n^{-1}d_n^{-2} S_{[nd_nt]} , t \in\,[0,1] \bigr\}
\]
converges in $%
D[0, 1]$ to $\sqrt{\eta}W_{H}$, where $\eta$ is a positive random
variable and $W_{H}$ is a standard fractional Brownian motion (with Hurst
index $H$), independent of $\eta$. This holds by item 1 of Lemma \ref
{AH} and the fact that $d_{n}^{-2}n^{-1}m_{[nd_{n}]}\rightarrow1$.
This completes the proof of Lemma \ref{AH}.
\end{pf}

We go back to the proofs of Theorems \ref{thmiso1}--\ref{thmiso3}. Note that the conditions of items 1 and 2 are
clearly satisfied by using either Comment \ref{SM} or \ref{comH} (with $
m_{n}=\sqrt{n}$), either Theorem \ref{IPlinearDMV} (with $m_{n}=\sqrt{n}
|\sum_{i=-n}^{n}a_{i}|$) or Theorem \ref{IP} or \ref{IP2} (with
$m_{n}=v_{n}$%
). In addition, in all these situations, we have that $m_{n}=(n^{\beta
}h(n))^{1/2}$ and the selection of $d_{n}$ leads to
\begin{eqnarray*}
d_{n}^{-2}n^{-1}m_{[nd_{n}]} &\sim&d_{n}^{(\beta-4)/2}n^{(\beta-2)/2}
\sqrt{h(nd_{n})} \\
&\sim& (L^{\ast}(n) )^{-1}\sqrt{h \bigl((nL^{\ast}(n))^{2/(4-\beta)} \bigr)} \\
&\sim& (L^{\ast}(n) )^{-1} (L(nL^{\ast}(n)) )^{-1} ,
\end{eqnarray*}
which converges to $1$ by (\ref{defL*}).
\end{pf*}
%s5 ###
\section{Proofs}

%s5.1 ###
\subsection{\texorpdfstring{Proof of Proposition \protect\ref{exp2}}{Proof of Proposition 2.1}}

Without loss of generality, we shall assume that $D_{\Psi}=1$ and $%
\sum_{j\in{\mathbf{Z}}}c_{m,j}^{2} =1$ since, otherwise, we can divide each
coefficient $c_{m,j}$ by $(\sum_{j\in{\mathbf{Z}}}c_{m,j}^{2})^{1/2}$ and
each variable by $D_{\Psi}$. Start with the decomposition
\begin{eqnarray*}
Y_{k}=\sum_{j=-\infty}^{\infty}P_{k-j}(Y_{k})=\sum_{j=-\infty
}^{\infty
}p_{j}P_{k-j}(Y_{k})/p_{j}.
\end{eqnarray*}
Then
\begin{eqnarray*}
S_{m}=\sum_{j=-\infty}^{\infty}p_{j}\sum_{k\in{\mathbf{Z}}%
}c_{m,k}P_{k-j}(Y_{k})/p_{j}.
\end{eqnarray*}

By using the facts that $\Psi$ is convex and non-decreasing, and
$p_{j}\geq
0$ with $\sum_{j\in{\mathbf{Z}}}p_{j}=D_{\Psi}=1$, we obtain that
\begin{eqnarray*}
\mathbf{E}\Psi(|S_{m}|)\leq\sum_{j=-\infty}^{\infty}p_{j}\mathbf
{E}\Psi
\biggl(\bigg|\sum_{k\in{\mathbf{Z}}}c_{m,k}P_{k-j}(Y_{k})/p_{j}\bigg| \biggr).
\end{eqnarray*}
Consider the martingale difference $U_{k}=c_{m,k}P_{k-j}(Y_{k})/p_{j}
$, $%
k\in{\mathbf{Z}}$. By Burkholder's inequality (see
%de la Pe\~{n}a and Gin\'{e}
\cite{DelapenaGine1999}, Theorem 6.6.2), we
obtain that
\begin{eqnarray*}
\mathbf{E}\Psi\biggl(\bigg|\sum_{k\in{\mathbf{Z}}}c_{m,k}P_{k-j}(Y_{k})/p_{j}\bigg|\biggr)\leq
K_{\alpha}\mathbf{E}\Psi\biggl(\biggl(\sum_{k\in{\mathbf{Z}}%
}c_{m,k}^{2}P_{k-j}^{2}(Y_{k})/p_{j}^{2}\biggr)^{1/2}\biggr) ,
\end{eqnarray*}
where $K_{\alpha}$ is a constant depending only on $\alpha$. Let $\Phi
(x)=\Psi(\sqrt{x})$. Since $\Phi$ is convex and $\sum_{k\in{\mathbf
{Z}}%
}c_{m,k}^{2}=1$, it follows that
\begin{eqnarray*}
\mathbf{E}\Psi\biggl(\bigg|\sum_{k\in{\mathbf{Z}}}c_{m,k}P_{k-j}(Y_{k})/p_{j}\bigg|\biggr)
&\leq &K_{\alpha}\mathbf{E}\Phi\biggl(\sum_{k\in{\mathbf{Z}}%
}c_{m,k}^{2}P_{k-j}^{2}(Y_{k})/p_{j}^{2}\biggr) \\
&\leq& K_{\alpha}\sum_{k\in{\mathbf{Z}}}c_{m,k}^{2}\mathbf{E}\Phi
\bigl(P_{k-j}^{2}(Y_{k})/p_{j}^{2}\bigr) \\
&\leq& K_{\alpha}\sum_{k\in{\mathbf{Z}}}c_{m,k}^{2}\mathbf{E}\bigl(\Psi\bigl(|P_{k-j}(Y_{k})|/p_{j}\bigr)\bigr).
\end{eqnarray*}
Therefore,
\begin{eqnarray*}
\mathbf{E}\Psi(|S_{m}|)\leq K_{\alpha}\sum_{k\in{\mathbf{Z}}%
}c_{m,k}^{2}\sum_{j=-\infty}^{\infty}p_{j} \mathbf{E}\bigl(\Psi
\bigl(|P_{k-j}(Y_{k})|/p_{j}\bigr)\bigr).
\end{eqnarray*}
Now, note that $\|P_{k-j}(Y_{k})\|_{\Psi}\leq p_{j}$, so using
the fact that $\sum_{k\in{\mathbf{Z}}}c_{m,k}^{2}=1$ and
$D_{\Psi}=\break\sum_{j=-\infty}^{\infty}p_{j}=1$, we get
\begin{eqnarray*}
\mathbf{E}\Psi(|S_{m}|)\leq K_{\alpha}
\end{eqnarray*}
and hence the desired result.

%s5.2 ###
\subsection{\texorpdfstring{Proof of Proposition \protect\ref{exp3}}{Proof of Proposition 2.2}}

Fix a positive integer $m$ and define
\begin{eqnarray*}
\theta_{0,m}=\sum_{k=0}^{2m-2}\sum_{i=k-m+1}^{m-1}P_{i}(\xi_{k}) \quad\mbox{and}\quad \theta_{j,m}=\theta_{0,m}\circ T^{j}.
\end{eqnarray*}
Observe that, by stationarity,
\begin{eqnarray*}
\Vert\theta_{0,m}\Vert_{\Psi}=\Bigg\Vert
\sum_{k=0}^{2m-2}\sum_{i=k-m+1}^{m-1}P_{i}(\xi_{k})\Bigg\Vert_{\Psi}\leq
2m\sum_{i\in{\mathbf{Z}}}\Vert P_{0}(\xi_{i})\Vert_{\Psi}<\infty.
\end{eqnarray*}
Simple computations lead to the decomposition
\begin{eqnarray*}
\sum_{i=-m+1}^{m-1}P_{i}(\xi_{0})-\sum_{\ell=1}^{2m-1}P_{m}(\xi_{\ell
})=\theta_{0,m}-\theta_{1,m} ,
\end{eqnarray*}
implying that%
\begin{eqnarray*}
\xi_{0}-\biggl(\sum_{k}P_{0}(\xi_{k})\biggr)\circ T^{m}=\theta_{0,m}-\theta
_{1,m}+\sum_{|i|\geq m}P_{i}(\xi_{0})-\biggl(\sum_{|k|\geq m}P_{0}(\xi
_{k})\biggr)\circ T^{m}\mbox{.}
\end{eqnarray*}
With our notation $(d_{0}=\sum_{k}P_{0}(\xi_{k}))$, we obtain%
%e29 ###
%
\begin{equation}\label{martdec}
\xi_{0}-d_{0}=d_{0}\circ T^{m}-d_{0}+\theta_{0,m}-\theta
_{1,m}+\sum_{|i|\geq m}P_{i}(\xi_{0})-\biggl(\sum_{|k|\geq m}P_{0}(\xi
_{k})\biggr)\circ T^{m}.
\end{equation}
By stationarity, we obtain similar decompositions for each $\xi_{j}-d_{j}.$
We shall treat the terms from the error of approximation $\sum_{i\in{%
\mathbf{Z}}}c_{n,i}(\xi_{i}-d_{i})$ separately. First, note that%
\begin{eqnarray*}
R_{1}:=\sum_{j=-{\infty}}^{\infty}c_{n,j}( d_{j}\circ
T^{m}-d_{j})&=&\sum_{j=-{\infty}}^{\infty
}(c_{n,j-m}-c_{n,j})d_{j}\\
&=&\sum_{k=0}^{m-1}\sum_{j=-{\infty}}^{\infty
}(c_{n,j-k-1}-c_{n,j-k})d_{j}.
\end{eqnarray*}
According to Proposition \ref{exp2},
\begin{eqnarray*}
\Vert R_{1}\Vert_{\Psi}\leq C_{\alpha}m\Vert d_{0}\Vert_{\Psi}
\Biggl(\sum_{j=-{\infty}}^{\infty}(c_{n,j}-c_{n,j-1})^{2}\Biggr)^{1/2}.
\end{eqnarray*}
To treat the second difference in the error, note that
\begin{eqnarray*}
R_{2}:=\sum_{i=-{\infty}}^{\infty}c_{n,i}(\theta_{i,m}-\theta
_{i+1,m})=\sum_{i=-{\infty}}^{\infty}(c_{n,i}-c_{n,i-1})\theta_{i,m}.
\end{eqnarray*}
By the definition of $\theta_{0,m}$, we have that
\begin{eqnarray*}
\sum_{j\in{\mathbf{Z}}}\Vert P_{j}(\theta_{0,m})\Vert_{\Psi}\leq
\sum_{k=0}^{2m-2}\sum_{i=k-m+1}^{m-1}\sum_{j\in{\mathbf{Z}}}\Vert
P_{j}(P_{i}(\xi_{k}))\Vert_{\Psi}.
\end{eqnarray*}
Now, $P_{j}(P_{i}(f))=0$ for $j\neq i$. It follows that
\begin{eqnarray*}
\sum_{j\in{\mathbf{Z}}}\Vert P_{j}(\theta_{0,m})\Vert_{\Psi}\leq
\sum_{k=0}^{2m-2}\sum_{\ell=k-m+1}^{m-1}\Vert P_{0}(\xi_{\ell})\Vert
_{\Psi}\leq(2m-1)\sum_{\ell=-m+1}^{m-1}\Vert P_{0}(\xi_{\ell})\Vert
_{\Psi}
\end{eqnarray*}
and, by Proposition \ref{exp2}, we conclude that
\begin{eqnarray*}
\|R_{2}\|_{\Psi}\leq2C_{\alpha}m\Biggl(\sum_{j=-{\infty}}^{\infty
}(c_{n,j}-c_{n,j-1})^{2}\Biggr)^{1/2}\sum_{\ell\in{\mathbf{Z}}}\Vert
P_{0}(\xi
_{\ell})\Vert_{\Psi}.
\end{eqnarray*}
For the term $R_{3}:=\sum_{i=-{\infty}}^{\infty}c_{n,i}(\sum_{|j|\geq
m}P_{j}(\xi_{0}))\circ T^{i}$, we apply Proposition \ref{exp2} to get%
\begin{eqnarray*}
\|R_{3}\|_{\Psi}\leq C_{\alpha}\Biggl(\sum_{i=-{\infty}}^{\infty
}c_{n,i}^{2}\Biggr)^{1/2} \sum_{|j|\geq m}\Vert P_{j}(\xi_{0})\Vert_{\Psi}.
\end{eqnarray*}
To deal with the last term $R_{4}:=\sum_{i=-{\infty}}^{\infty
}c_{n,i}(\sum_{|k|\geq m}P_{0}(\xi_{k}))\circ T^{m+i}$, we again apply
Proposition~\ref{exp2}, which gives
\begin{eqnarray*}
\|R_{4}\|_{\Psi}\leq C_{\alpha}\Biggl(\sum_{i=-{\infty}}^{\infty
}c_{n,i}^{2}\Biggr)^{1/2}\sum_{|k|\geq m}\|P_{0}(\xi_{k})\Vert_{\Psi}.
\end{eqnarray*}
Combining all the bounds, we obtain the desired approximation.

%s5.3 ###
\subsection{\texorpdfstring{Proof of Lemma \protect\ref{max}}{Proof of Lemma 2.1}}

For any $m\in[1,2^{N}]$, write $m$ in base $2$ as follows:
\begin{eqnarray*}
m=\sum_{i=0}^{N}b_{i}(m)2^{i},\qquad\mbox{where $b_{i}(m)=0$ or $b_{i}(m)=1$}.
\end{eqnarray*}
Set $m_{L}=\sum_{i=L}^{N}b_{i}(m)2^{i}$. So, for any $p\geq1$, we have
\begin{eqnarray*}
|S_{m}|^{p}\leq \Biggl(\sum_{L=0}^{N}|S_{m_{L}}-S_{m_{L+1}}| \Biggr)^{p}.
\end{eqnarray*}
Hence, setting
\begin{eqnarray*}
\alpha_{L}=\|S_{2^{L}}\|_{\Psi_{p}} (\Psi^{-1}(2^{N-L}) )^{1/p}%
\quad\mbox{and}\quad\lambda_{L}=\frac{\alpha_{L}}{\sum_{L=0}^{N}\alpha_{L}} ,
\end{eqnarray*}
we get, by convexity,
\begin{eqnarray*}
|S_{m}|^{p}\leq\sum_{L=0}^{N}\lambda
_{L}^{1-p}|S_{m_{L}}-S_{m_{L+1}}|^{p}.
\end{eqnarray*}
Now, $m_{L}\neq m_{L+1}$ only if $b_{L}(m)=1$ and, in that case, $%
m_{L}=k_{m}2^{L}$ with $k_{m}$ odd. It follows that
\begin{eqnarray*}
\max_{1\leq m\leq2^{N}}|S_{m}|^{p}\leq\sum_{L=0}^{N}\lambda
_{L}^{1-p}\max_{1\leq k\leq2^{N-L},k\ \mathrm{odd}}\big|S_{k2^{L}}-S_{(k-1)2^{L}}\big|^{p}.
\end{eqnarray*}
Now, we apply \cite{LedouxTalagrand1991},
Lemma 11.3, to the
variables
\begin{eqnarray*}
Z_{k}=\frac{|S_{k2^{L}}-S_{(k-1)2^{L}}|^{p}}{A^{p}},\qquad\mbox{ where $A=\|S_{2^{L}}\|_{\Psi_{p}}$,}
\end{eqnarray*}
and to the Young function $\Psi$. Since
\begin{eqnarray*}
{\mathbf{E}}(\Psi(Z_{k}))={\mathbf{E}}\Psi_{p} \biggl(\frac{|S_{2^{L}}|}{A}%
\biggr)\leq1
\end{eqnarray*}
and since $\Psi^{-1}$ is concave, we get that, for any measurable
set $B$,
\begin{eqnarray*}
{\mathbf{E}}(Z_{k}{\mathbf{1}}_{B})\leq P(B)\Psi^{-1} \biggl(\frac{1}{P(B)}%
\biggr)
\end{eqnarray*}
so that the assumptions of Ledoux and Talagrand \cite{LedouxTalagrand1991}, Lemma 11.3, are
satisfied. It follows that
\begin{eqnarray*}
{\mathbf{E}} \Bigl(\max_{1\leq k\leq2^{N-L},k\ \mathrm{odd}}\big|S_{k2^{L}}-S_{(k-1)2^{L}}\big|^{p} \Bigr)\leq A^{p}\Psi^{-1}(2^{N-L}).
\end{eqnarray*}
Finally, we conclude that
\begin{eqnarray*}
{\mathbf{E}}\Bigl(\max_{1\leq m\leq2^{N}}|S_{m}|^{p}\Bigr)\leq \Biggl(%
\sum_{L=0}^{N}\alpha_{L} \Biggr)^{p} ,
\end{eqnarray*}
which is the desired result.

%s5.4 ###
\subsection{\texorpdfstring{Proofs of Theorems \protect\ref{IP} and
\protect\ref{IP2}}{Proofs of Theorems 3.1 and 3.2}}

By the weak convergence theory of random functions, it suffices to establish
the convergence of the finite-dimensional distributions and the
tightness of
$\{v_n^{-1}S_{[nt]}, t\in[0,1]\}$. For the finite-dimensional
distribution, we shall use the following proposition which was basically
established in \cite{PeligradUtev1997,PeligradUtev2006b}.
%p5.1
\begin{proposition}
\label{propCLTtri} Let $\{\xi_{k}\}_{k\in\mathbf{Z}}$ be a strictly
stationary sequence of centered and regular random variables in
${\mathbf{L}}%
^{2}$ such that $\sum_{j}\Vert P_{0}(\xi_{j})\Vert_{2}<\infty$. For any
positive integer $n$, let $\{b_{n,i},-\infty\leq i\leq\infty\}$ be a
triangular array of numbers satisfying
%e30 ###
%
\begin{equation}\label{cond1cni}
\sum_{i}b_{n,i}^{2}\rightarrow1 \quad\mbox{and}\quad \sum
_{j}(b_{n,j}-b_{n,j-1})^{2}\rightarrow0 \qquad\mbox{as }n\rightarrow\infty
\end{equation}
and
%e31 ###
%
\begin{equation}\label{cond2cni**}
\sup_{j}|b_{n,j}|\rightarrow0 \qquad\mbox{as } n\rightarrow
\infty.
\end{equation}
Then $\{S_{n}=\sum_{j}b_{n,j}\xi_{j}\}$ converges in distribution to
$\sqrt{%
\eta}N$, where $N$ is a standard Gaussian random variable, independent
of $%
\eta$, and $\eta=\sum_{k \in{\mathbf{Z}}}{\mathbf{E}}(\xi_{0}\xi
_{k}|{%
\mathcal{I}})$.
\end{proposition}

\begin{pf}
We give here the
proof for completeness. By using Proposition \ref{exp3}, it suffices to prove
this proposition with $d_{j}=d_{0}\circ T^{j}$ in place of $\xi_{j}$, where
$d_{0}=\sum_{j}P_{0}(\xi_{j})$. Hence, we just have to apply the central
limit theorem for triangular arrays of martingales (see \cite{HallHeyde1980}, Theorem 3.6). The Lindeberg
condition has been established by Peligrad
and Utev \cite{PeligradUtev1997}, provided that condition~(\ref
{cond2cni**}) and the first part of
condition (\ref{cond1cni}) are satisfied. Now, in the proof of Peligrad
and Utev \cite{PeligradUtev2006b},
Proposition 4, it is established that (\ref%
{cond1cni}) implies that
\begin{eqnarray*}
\sum_{j}b_{n,j}^{2}d_{j}^{2}\rightarrow\eta\qquad \mbox{in probability as $
n\rightarrow\infty$,}
\end{eqnarray*}
which ends the proof of the proposition.
\end{pf}

We return to the proofs of Theorems \ref{IP} and \ref{IP2}. To prove the
convergence of the finite-dimensional distributions, we shall apply the
Cram%
\'{e}r--Wold device. For all integer $1\leq\ell\leq m$, let $n_{\ell
}=[nt_{\ell}]$, where $0<t_{1}<t_{2}<\cdots<t_{m}\leq1$. For $\lambda
_{1},\ldots,\lambda_{m}\in{\mathbf{R}}$, note that
%e32 ###
%
\begin{equation}
\frac{\sum_{\ell=1}^{m}\lambda_{\ell}S_{n_{\ell}}}{v_{n}}=\sum
_{j\in{%
\mathbf{Z}}} \Biggl(\sum_{\ell=1}^{m}\frac{\lambda_{\ell}c_{n_{\ell},j}}{%
v_{n}} \Biggr)\xi_{j} , \label{decCW}
\end{equation}
where $c_{n,j}=a_{1-j}+\cdots+a_{n-j}$ for all $j\in{\mathbf{Z}}$ and
$%
v_{n}^{2}=\sum_{j\in{\mathbf{Z}}}c_{n,j}^{2}$. Let
%e33 ###
%
\begin{equation}
b_{n,j}=\frac{1}{\Lambda_{m,\beta}}\sum_{\ell=1}^{m}\frac{\lambda
_{\ell
}c_{n_{\ell},j}}{v_{n}} , \label{defcnj}
\end{equation}
where
\begin{eqnarray*}
\Lambda_{m,\beta}^{2}=\frac{1}{2}\sum_{\ell,k=1}^{m}\lambda_{\ell
}\lambda_{k} (t_{\ell}^{\beta}+t_{k}^{\beta}-|t_{k}-t_{\ell
}|^{\beta} ).
\end{eqnarray*}
We apply Proposition \ref{propCLTtri} to $b_{n,j}$ and the $\xi_{j}$'s
defined as $\Lambda_{m,\beta}\xi_{j}$. First, we have to calculate the
limit over $n$ of the quantity
\begin{eqnarray*}
\sum_{j\in{\mathbf{Z}}}b_{n,j}^{2}=\frac{1}{\Lambda_{m,\beta
}^{2}}\frac{%
\sum_{j\in{\mathbf{Z}}}\sum_{\ell=1}^{m}\sum_{k=1}^{m}\lambda_{\ell
}\lambda_{k}c_{n_{\ell},j}c_{n_{k},j}}{v_{n}^{2}}.
\end{eqnarray*}
For any $1\leq{\ell}\leq k\leq m$, by using the fact that for any two real
numbers $A$ and $B$, we have $A(A+B)=1/2(A^{2}+(A+B)^{2}-B^{2})$, we
get that
\begin{eqnarray*}
\frac{1}{v_{n}^{2}}\sum_{j\in{\mathbf{Z}}}c_{n_{\ell},j}c_{n_{k},j}
&=&%
\frac{1}{2v_{n}^{2}}\sum_{j\in{\mathbf{Z}}} \bigl(c_{n_{\ell
},j}^{2}+c_{n_{k},j}^{2}- (c_{n_{\ell},j}-c_{n_{k},j} )^{2} \bigr) \\
&=&\frac{1}{2v_{n}^{2}}\sum_{j\in{\mathbf{Z}}} (c_{n_{\ell
},j}^{2}+c_{n_{k},j}^{2}-c_{n_{k}-n_{\ell},j}^{2} ).
\end{eqnarray*}
By now using condition (\ref{hyposn}), we derive that, for any $1\leq{%
\ell}\leq k\leq m$,
%e34 ###
%
\begin{equation}
\frac{\sum_{j\in{\mathbf{Z}}}b_{n_{\ell},j}b_{n_{k},j}}{v_{n}^{2}}%
\rightarrow\frac{1}{2} \bigl(t_{\ell}^{\beta}+t_{k}^{\beta
}-(t_{k}-t_{\ell})^{\beta} \bigr). \label{doubleprod}
\end{equation}
It follows from (\ref{doubleprod}) that
%e35 ###
%
\begin{equation}
\lim_{n\rightarrow\infty}\sum_{j\in{\mathbf{Z}}}b_{n,j}^{2}=1.
\label{limcnj}
\end{equation}
As a consequence, the first part of condition (\ref{cond1cni}) holds.
On the
other hand, by using Peligrad and Utev \cite{PeligradUtev2006b},
Lemma A.1, the second
part of condition (\ref{cond1cni}) is satisfied. Now, by the proof of
Corollary 2.1 in \cite{PeligradUtev1997}, we
get that
\begin{eqnarray*}
\frac{\max_{j}|c_{n,j}|}{v_{n}}\rightarrow0 ,
\end{eqnarray*}
which, together with (\ref{hyposn}), implies (\ref{cond2cni**}). Now, applying
Proposition \ref{propCLTtri}, we derive that
\begin{eqnarray*}
\frac{\sum_{\ell=1}^{m}\lambda_{\ell}S_{n_{\ell}}}{v_{n}}\qquad\mbox{converges in distribution to }\Lambda_{m,\beta}\sqrt{\eta}N ,
\end{eqnarray*}
ending the proof of the convergence of the finite-dimensional distribution.

We now turn to the proof of the tightness of $\{v_{n}^{-1}S_{[nt]},t\in
[0,1]\}$. By using Proposition \ref{exp2}, we get, for $q\geq
2$, that
%e36 ###
%
\begin{equation}
\|S_{k}\|_{q}\leq C_{q} \biggl(\sum_{j\in{\mathbf{Z}}}b_{k,j}^{2} \biggr)%
^{1/2}\sum_{m\in{\mathbf{Z}}}\Vert P_{0}(\xi_{m})\Vert
_{q}=C_{q}v_{k}\sum_{m\in{\mathbf{Z}}}\Vert P_{0}(\xi_{m})\Vert_{q} ,
\label{r1lemmamom}
\end{equation}
provided that $\sum_{m\in{\mathbf{Z}}}\Vert P_{0}(\xi_{m})\Vert
_{q}<\infty$. Therefore, the conditions of Taqqu \cite{Taqqu1975}, Lemma 2.1, page
290,
are satisfied with $q>2/\beta$ and the tightness follows.

Finally, to prove (\ref{inemawxIP}), we use (\ref{r1lemmamom}),
together with
Lemma \ref{max} applied with $\psi(x)=x$, by taking into account the
fact that $v_n^2$
is regularly varying with exponent $\beta$.

%s5.5 ###
\subsection{\texorpdfstring{Proof of Remarks \protect\ref{use} and
\protect\ref{corIPlinearDMV}}{Proof of Remarks 3.3 and 3.6}}

To prove Remark \ref{use}, we apply Lemma \ref{lnumbers} from the \hyperref[Appendix]{Appendix}
with $b_i=1$ and $u_{i}=\Vert P_{-i}(\xi_{0})\Vert_{q}$. Hence, we get
\begin{eqnarray*}
\sum_{n=1}^{\infty}\Vert P_{-n}(\xi_{0})\Vert_{q}\leq C_{q}\sum
_{n=1}^{\infty} \Biggl(\frac{1}{n}\sum_{k=n}^{\infty}\Vert P_{-k}(\xi
_{0})\Vert_{q}^{q} \Biggr)^{1/q}.
\end{eqnarray*}
Applying the Rosenthal inequality given in
\cite{HallHeyde1980}, Theorem 2.12, we then derive that for any
$q\in[2,\infty[$, there
exists a constant $c_{q}$, depending only on $q$, such that
\begin{eqnarray*}
\sum_{k=n}^{\infty}\Vert P_{-k}(\xi_{0})\Vert_{q}^{q}\leq c_{q}\Bigg\Vert
\sum_{k=n}^{\infty}P_{-k}(\xi_{0})\Bigg\Vert_{q}^{q}=c_{q}\Vert{\mathbf
{E}}%
(\xi_{n}|{\mathcal{F}}_{0})\Vert_{q}^{q}.
\end{eqnarray*}
The same argument works with $P_{-i}(\xi_{0})$ replaced by $P_{i}(\xi
_{0}) $, and the result follows by applying the Rosenthal inequality and
noting that $\Vert\xi_{-n}-{\mathbf{E}}(\xi_{-n}|{\mathcal{F}}%
_{0})\Vert_{q}=\Vert\sum_{k=n}^{\infty}P_{k+1}(\xi_{0})\Vert_{q}$.

To prove Remark \ref{corIPlinearDMV}, we apply Lemma \ref{lnumbers}
from the
\hyperref[Appendix]{Appendix} with $b_{n}=\log(n)$ and $u_{n}=\Vert P_{0}(\xi_{n})\Vert_{2}$.
We then get that
\begin{eqnarray*}
\sum_{n=1}^{\infty}\log n\Vert P_{0}(\xi_{n})\Vert_{2}\leq C
\sum_{n=1}^{\infty}\frac{\log n}{\sqrt{n}} \Biggl(\sum_{k=n}^{\infty}\Vert
P_{0}(\xi_{k})\Vert_{2}^{2} \Biggr)^{1/2}.
\end{eqnarray*}
Now, note that
\begin{eqnarray*}
\sum_{k=n}^{\infty}\Vert P_{0}(\xi_{k})\Vert_{2}^{2}=\Vert{\mathbf
{E}}%
(\xi_{n}|{\mathcal{F}}_{0})\Vert_{2}^{2}
\end{eqnarray*}
and so
\begin{eqnarray*}
\sum_{n=1}^{\infty}\log n\Vert P_{0}(\xi_{n})\Vert_{2}\leq
C\sum_{n=1}^{\infty}\log n\frac{\Vert{\mathbf{E}}(\xi_{n}|\mathcal{F}
_{0})\Vert_{2}}{\sqrt{n}}<\infty.
\end{eqnarray*}
The same argument works with $P_{0}(\xi_{i})$ replaced by $P_{0}(\xi
_{-i})$%
.

%s5.6 ###
\subsection{\texorpdfstring{Proof of Theorem
\protect\ref{IPlinearDMV}}{Proof of Theorem 3.3}}

For all $j\in{\mathbf{Z}}$, let $d_{j}=\sum_{\ell\in{\mathbf{Z}}%
}P_{j}(\xi_{\ell})$. Note that if either condition (a) or
condition (b) is satisfied, $(d_{j})_{j\in{\mathbf{Z}}}$ is a
sequence of martingale differences in ${\mathbf{L}}^{2}$. We set
\begin{eqnarray*}
Y_{k}=\sum_{i\in{\mathbf{Z}}}a_{i}d_{k-i}\quad\mbox{and}\quad T_{n}=%
\sum_{k=1}^{n}Y_{k} ,
\end{eqnarray*}
and apply \cite{DMPU},
Corollary 4. By
taking into account Remark \ref{rmsn}, we derive that under (\ref{condDMV}),
\begin{eqnarray*}
\bigl\{s_{n}^{-1}T_{[nt]},t\in[0,1]\bigr\}\quad\mbox{converges in distribution in
$(D([0,1]),d)$ to \quad$\sqrt{{\mathbf{E}} (d_{0}^{2}|{\mathcal{I}} )} W $,}
\end{eqnarray*}
where $W$ is a standard Brownian motion independent of ${\mathcal{I}}$. It
follows that in order to prove that $\{s_{n}^{-1}S_{[nt]},t\in[
0,1]\} $ converges in distribution in $(D([0,1]),d)$ to $\sqrt{{\mathbf
{E}}%
(d_{0}^{2}|{\mathcal{I}} )}W$, it is sufficient to show that
%e37 ###
%
\begin{equation}
\frac{\Vert\max_{1\leq k\leq n}|S_{k}-T_{k}|\Vert
_{2}}{s_{n}}\rightarrow0%
\qquad\mbox{as $n\rightarrow\infty$}. \label{approxIP*}
\end{equation}
Now, for any $n$, let $N$ be such that $2^{N-1}<n\leq2^{N}$. By using
Remark %
\ref{rmsn} and the properties of the slowly varying function, we get
that $%
s_{n}\sim s_{2^{N}}$. So, the proof (\ref{approxIP*}) is reduced to showing
that
%e38 ###
%
\begin{equation}
\frac{\Vert\max_{1\leq k\leq2^{N}}|S_{k}-T_{k}|\Vert_{2}}{s_{2^{N}}}%
\rightarrow0\qquad\mbox{as $N\rightarrow\infty$}. \label{approxIP}
\end{equation}
We first prove that (\ref{approxIP}) holds under condition (a). By using
Corollary \ref{corexp3}, together with Lemma~\ref{max}, we get that for any
positive integer $m$,
\begin{eqnarray*}
\Big\Vert\max_{1\leq k\leq2^{N}}|S_{k}-T_{k}|\Big\Vert_{2}
&\leq&C_{1}\sum_{|k|\geq m}\Vert P_{0}(\xi_{k})\Vert_{\Psi_{2,\alpha
}}\sum_{L=0}^{N}v_{2^{L}} (g^{-1}(2^{N-L}) )^{1/2} \\
&&{} +C_{2}m\sum_{L=0}^{N} (g^{-1}(2^{N-L}) )^{1/2} ,
\end{eqnarray*}
where $g(x)=x\log^{\alpha}(1+x)$. Noting that $g^{-1}(x)\sim\frac{x%
}{\log^{\alpha}(1+x)}$ as $x$ goes to infinity, and taking into account
Remark \ref{rmsn} and the first part of condition (\ref{condDMV}), we get
that
%e39 ###
%
\begin{equation}\label{ineapproxIP}
\Big\Vert\max_{1\leq k\leq2^{N}}|S_{k}-T_{k}|\Big\Vert_{2}\leq
Cs_{2^{N}}\sum_{|k|\geq m}\Vert P_{0}(\xi_{k})\Vert_{\Psi_{2,\alpha
}}+Cm\epsilon(N)s_{2^{N}} ,
\end{equation}
where $\epsilon(N)\rightarrow0$ as $N\rightarrow\infty$. By now
using (%
\ref{ineapproxIP}) and first letting $N$ tend to infinity and then $m$ tend
to infinity, we derive (\ref{approxIP}) under condition (a).

We now turn to the proof of (\ref{approxIP}) under condition (b).
Taking $%
m=m_{2^{L}}=2^{L/4}$ in Corollary \ref{corexp3} and using Lemma \ref{max}
with $p=2$ and $\psi(x)=x$, we get that
\begin{eqnarray}\label{approx1Cb}
&&\frac{\Vert\max_{1\leq k\leq2^{N}}|S_{k}-T_{k}|\Vert_{2}}{s_{2^{N}}}\nonumber
\\[-8pt]\\[-8pt]
&&\quad\leq C
\frac{2^{N/2}}{s_{2^{N}}}\sum_{L=0}^{N}\frac{m_{2^{L}}}{2^{L/2}}\nonumber
+C\frac{2^{N/2}}{s_{2^{N}}}\sum_{L=0}^{N}\frac{v_{2^{L}}}{2^{L/2}}%
\sum_{|k|\geq m_{2^{L}}}\Vert P_{0}(\xi_{k})\Vert_{2}.\nonumber
\end{eqnarray}
By Remark \ref{rmsn}, we have that $\lim_{N\rightarrow\infty}\frac
{s_{2^{N}}%
}{2^{N/2}}=\infty$, which, together with the selection of $m_{2^{L}}$, implies
that the first term on the right-hand side of the above inequality
tends to zero
as $n\rightarrow\infty$. Now, to treat the last term, we first fix a
positive integer $p $ and write
\begin{eqnarray*}
\frac{2^{N/2}}{s_{2^{N}}}\sum_{L=0}^{N}\frac{v_{2^{L}}}{2^{L/2}}%
\sum_{|k|\geq m_{2^{L}}}\Vert P_{0}(\xi_{k})\Vert_{2} &\leq& p \frac{
2^{N/2}}{s_{2^{N}}}\max_{0\leq L < p }\frac{v_{2^{L}}}{2^{L/2}}\sum
_{|k|\geq
m_{2^{L}}}\Vert P_{0}(\xi_{k})\Vert_{2} \\
&&{} +\frac{2^{N/2}}{s_{2^{N}}}\sum_{L=p }^{N}\frac{v_{2^{L}}}{2^{L/2}}%
\sum_{|k|\geq m_{2^{L}}}\Vert P_{0}(\xi_{k})\Vert_{2}.
\end{eqnarray*}
Since $\lim_{N\rightarrow\infty}\frac{s_{2^{N}}}{2^{N/2}}=\infty$, the
first term on the right-hand side of the above inequality tends to zero
as $%
N\rightarrow\infty$. To treat the second one, we note that if
$N$ and $p$ are large enough,
\begin{eqnarray*}
\frac{2^{N/2}}{s_{2^{N}}}\sum_{L=p }^{N}\frac{v_{2^{L}}}{2^{L/2}}%
\sum_{|k|\geq m_{2^{L}}}\Vert P_{0}(\xi_{k})\Vert_{2}\leq C\sum_{L=p
}^{N}%
\frac{h(2^{L})}{h(2^{N})}\sum_{|k|\geq m_{2^{L}}}\Vert P_{0}(\xi
_{k})\Vert
_{2} ,
\end{eqnarray*}
where $h(n)= |\sum_{i=-n}^{n}a_{i} |$. By the first part of condition (%
\ref{condDMV}),
\begin{eqnarray*}
\limsup_{N\rightarrow\infty}\max_{p \leq L\leq N}\frac
{h(2^{L})}{h(2^{N})}%
<\infty.
\end{eqnarray*}
Hence, for $N$ and $p $ large enough and taking into account the
selection of $m_{2^{L}}$, we get that
\begin{eqnarray*}
\frac{2^{N/2}}{s_{2^{N}}}\sum_{L=p }^{N}\frac{v_{2^{L}}}{2^{L/2}}%
\sum_{|k|\geq m_{2^{L}}}\Vert P_{0}(\xi_{k})\Vert_{2}\leq C\sum
_{|k|\geq{%
2^{p /4}}}\log k\Vert P_{0}(\xi_{k})\Vert_{2} ,
\end{eqnarray*}
which converges to zero as $p \rightarrow\infty$, by using condition (b).
Hence, starting from (\ref{approx1Cb}) and taking into account the previous
considerations, we get that (\ref{approxIP}) holds under condition (b). The
proof of (\ref{inemawxIPlinearDMV}) is straightforward, following the
arguments used
to derive (\ref{approxIP*}).

%s5.7 ###
\subsection{\texorpdfstring{Proof of Comment \protect\ref{comH}}{Proof of Comment 3.2}}

The justification of this result is due to the following coboundary
decomposition. Define
%e40 ###
%
\begin{equation}
Z_{0}=\sum_{\ell=1}^{\infty}\sum_{k=\ell}^{\infty}a_{k}\xi_{-\ell
}-\sum_{\ell=0}^{\infty}\sum_{k=-\infty}^{-\ell-1}a_{k}\xi_{\ell}.
\label{defcobord}
\end{equation}
Since condition (\ref{HeyHan}) implies that the sequence $(\xi
_{i})_{i\in{%
\mathbf{Z}}}$ has a bounded spectral density, the random variable
$Z_{0}$ is well defined in ${\mathbf{L}}^{2}$ under condition
(H). Now,
\begin{eqnarray*}
Z_{0}-Z_{0}\circ T=\sum_{\ell=1}^{\infty}a_{\ell}\xi_{-\ell}-\xi
_{0}\sum_{k=1}^{\infty}a_{k}-\xi_{0}\sum_{k=1}^{\infty}a_{-k}+\sum
_{\ell
=1}^{\infty}a_{-\ell}\xi_{\ell} ,
\end{eqnarray*}
whence
\begin{eqnarray*}
A\xi_{0}+Z_{0}-Z_{0}\circ T=a_{0}\xi_{0}+\sum_{j\in{\mathbf
{Z}}\setminus
\{0\}}a_{j}\xi_{-j}=X_{0}.
\end{eqnarray*}
We derive that, for any $k\geq1$,
%e41 ###
%
\begin{equation}
S_{k}=A\sum_{i=1}^{k}\xi_{i}+Z_{1}-Z_{k+1} , \label{Sncobord}
\end{equation}
where $Z_{k}=Z_{0}\circ T_{k}$. Since, under condition (\ref{HeyHan}),
the partial sums
process $\{n^{-1/2}\sum_{k=1}^{[nt]}\xi_{k},t\in[0,1]\}$ converges
in distribution in $D([0,1])$ to $\sqrt{\lambda}W$ with $\lambda
=\sum_{j\in{\mathbf{Z}}}{\mathbf{E}}(\xi_{0}\xi_{j}|{\mathcal{I}})$, we
just have to show that
\begin{eqnarray*}
\limsup_{n\rightarrow\infty}{\mathbf{P}} \Bigl(\max_{1\leq k\leq
n}|Z_{k+1}|\geq\varepsilon\sqrt{n} \Bigr)=0 ,
\end{eqnarray*}
which holds because $Z_{0}\in{\mathbf{L}}^{2}$ (see \cite{HallHeyde1980}, inequality
(5.30)).

\begin{appendix}
%s6 ###
\section*{Appendix}\label{Appendix}
\setcounter{subsection}{0}
%s6.1 ###
\subsection{A fact concerning series}
%lA.1
\begin{lemma}
\label{lnumbers} Let $q>1$ and $\alpha= 2(q-1)/q$. Let
$(b_{j})_{j\in{\mathbf{N}}}$ be a sequence of non-negative numbers
such that $n^{\alpha} b_n \leq
K_{\alpha} \sum_{k=1}^n k^{\alpha-1} b_k$ for some positive constant $%
K_\alpha$ depending only on $\alpha$. Then, for any sequence of
non-negative numbers $(u_{j})_{j\in{\mathbf{N}}}$, the following
inequality holds:
\begin{eqnarray*}
\sum_{n=1}^{\infty}b_n u_{n}\leq C_{q} \sum_{n=1}^{\infty}b_n \Biggl(
\frac{1}{n}\sum_{k=n}^{\infty}u_{k}^{q} \Biggr) ^{1/q},
\end{eqnarray*}
where $C_{q}$ is a constant depending only on $q$.
\end{lemma}

\begin{pf}
We write
\begin{eqnarray*}
\sum_{n=1}^{\infty}b_{n}u_{n} \leq K_\alpha\sum_{n=1}^{\infty
}n^{-\alpha}u_{n} \Biggl(\sum_{k=1}^{n} b_{k}k^{\alpha-1} \Biggr)\leq
K_\alpha\sum_{k=1}^{\infty}b_{k}k^{\alpha-1} \biggl(\sum_{n\geq k}
n^{-\alpha}u_{n} \biggr).
\end{eqnarray*}
H\"{o}lder's inequality then gives
\begin{eqnarray*}
\sum_{n=1}^{\infty}b_{n}u_{n}\leq C_{q}^{\prime}\sum_{k=1}^{\infty
}b_{k}k^{\alpha-1} \biggl(\sum_{n\geq k} n^{-2} \biggr)^{\alpha/2} \biggl(%
\sum_{n\geq k} u_{n}^{q} \biggr)^{1/q}
\end{eqnarray*}
and the result follows.
\end{pf}
\end{appendix}

\section*{Acknowledgements}
The authors would like to thank Sergey Utev for valuable discussions on this
topic. They are also grateful to the referees and the Associate Editor
for helpful comments that
improved the presentation of this paper.
The research of Magda Peligrad was supported in part by a Charles
Phelps Taft Memorial Fund grant and NSA Grant No.~H98230-09-1-0005.

\printhistory

\end{document}